\documentclass[final,leqno,onefignum,onetabnum]{siamltex1213}
\usepackage{amsfonts}
\usepackage{amsmath,booktabs,ctable,threeparttable}
\usepackage{amssymb,amsfonts,boxedminipage}
\newtheorem{remark}{Remark}[section]

\title{Regularization Properties of the Krylov Iterative Solvers
CGME and LSMR For Linear Discrete Ill-Posed
Problems with an Application to Truncated Randomized SVDs\thanks{This
work was supported in part by
the National Science Foundation of China (No. 11771249)}}

\author{Zhongxiao Jia\thanks{Department of Mathematical Sciences, Tsinghua
University, 100084 Beijing, China. (\email{jiazx@tsinghua.edu.cn})}}

\begin{document}
\maketitle
\slugger{sirev}{xxxx}{xx}{x}{x--x}

\begin{abstract}
For the large-scale linear discrete ill-posed problem $\min\|Ax-b\|$ or $Ax=b$
with $b$ contaminated by Gaussian white noise, there are four commonly used
Krylov solvers: LSQR and its mathematically equivalent CGLS,
the Conjugate Gradient (CG) method applied to $A^TAx=A^Tb$, CGME, the CG method
applied to $\min\|AA^Ty-b\|$ or $AA^Ty=b$ with $x=A^Ty$, and
LSMR, the minimal residual (MINRES) method
applied to $A^TAx=A^Tb$. These methods have intrinsic regularizing effects,
where the number $k$ of iterations plays the role of the regularization parameter.
In this paper, we establish a number of regularization properties
of CGME and LSMR, including the filtered SVD expansion of CGME iterates, and prove
that the 2-norm filtering best regularized solutions by CGME and LSMR
are less accurate than and at least as accurate as those by LSQR, respectively.
We also prove that the semi-convergence of CGME and LSMR always occurs no later
and sooner than that of LSQR, respectively.
As a byproduct, using the analysis approach for CGME,
we improve a fundamental result on the
accuracy of the truncated rank $k$ approximate SVD of $A$ generated by
randomized algorithms, and reveal how the truncation step damages the accuracy.
Numerical experiments justify our results on CGME and LSMR.
\end{abstract}



\pagestyle{myheadings}
\thispagestyle{plain}
\markboth{ZHONGXIAO JIA}{REGULARIZATION OF CGME AND LSMR}

\section{Introduction and Preliminaries}\label{intro}

Consider the linear discrete ill-posed problem
\begin{equation}
  \min\limits_{x\in \mathbb{R}^{n}}\|Ax-b\| \mbox{\,\ or \ $Ax=b$,}
  \ \ \ A\in \mathbb{R}^{m\times n}, \label{eq1}
  \ b\in \mathbb{R}^{m},
\end{equation}
where the norm $\|\cdot\|$ is the 2-norm of a vector or matrix, and
$A$ is extremely ill conditioned with its singular values decaying
to zero without a noticeable gap. We simply assume that $m\geq n$.
Since the results in this paper hold for both the $m\geq n$ and $m\leq n$
cases. \eqref{eq1} arises from many applications, e.g., from
the discretization of the first kind Fredholm integral equation
\begin{equation}\label{eq2}
Kx=(Kx)(t)=\int_{\Omega} k(s,t)x(t)dt=g(s)=g,\ s\in \Omega
\subset\mathbb{R}^q,
\end{equation}
where the kernel $k(s,t)\in L^2({\Omega\times\Omega})$ and
$g(s)$ are known functions, while $x(t)$ is the
unknown function to be sought.
Applications include image deblurring, signal processing, geophysics,
computerized tomography, heat propagation, biomedical and optical imaging,
groundwater modeling, and many others
\cite{aster,engl93,engl00,hansen10,kaipio,kern,kirsch,natterer,vogel02}.
The right-hand side $b=b_{true}+e$ is assumed to be
contaminated by a Gaussian white noise $e$, caused by measurement, modeling
or discretization errors, where $b_{true}$
is noise-free and $\|e\|<\|b_{true}\|$.
Because of the presence of noise $e$ and the extreme
ill-conditioning of $A$, the naive
solution $x_{naive}=A^{\dagger}b$ of \eqref{eq1} generally bears no
relation to the true solution $x_{true}=A^{\dagger}b_{true}$, where
$\dagger$ denotes the Moore-Penrose inverse of a matrix.
Therefore, we must use regularization to extract a
good approximation to $x_{true}$ as much as possible.

For a Gaussian white noise $e$,
throughout the paper, we always assume that $b_{true}$ satisfies the discrete
Picard condition $\|A^{\dagger}b_{true}\|\leq C$ with some constant $C$
for $\|A^{\dagger}\|$ arbitrarily
large \cite{aster,gazzola15,hansen90,hansen90b,hansen98,hansen10,kern}.
Without loss of generality,
assume that $Ax_{true}=b_{true}$. Then a dominating regularization approach is
to solve the problem
\begin{equation}\label{posed}
\min\limits_{x\in \mathbb{R}^{n}}\|Lx\| \ \ \mbox{subject to}\ \
\|Ax-b\|\leq \tau\|e\|
\end{equation}
with $\tau>1$ slightly \cite{hansen98,hansen10},
where $L$ is a regularization matrix and its suitable choice is based on
a-prior information on $x_{true}$.

In this paper, we are concerned with the case $L=I$ in \eqref{posed}, 
which corresponds to a 2-norm filtering regularization problem. Let
\begin{equation}\label{eqsvd}
  A=U\left(\begin{array}{c} \Sigma \\ \mathbf{0} \end{array}\right) V^{T}
\end{equation}
be the singular value decomposition (SVD) of $A$,
where $U = (u_1,u_2,\ldots,u_m)\in\mathbb{R}^{m\times m}$ and
$V = (v_1,v_2,\ldots,v_n)\in\mathbb{R}^{n\times n}$ are orthogonal,
$\Sigma = {\rm diag} (\sigma_1,\sigma_2,\ldots,\sigma_n)\in\mathbb{R}^{n\times n}$
with the singular values
$\sigma_1>\sigma_2 >\cdots >\sigma_n>0$ assumed to be simple,
the superscript $T$
denotes the transpose of a matrix or vector, and
$\mathbf{0}$ denotes a zero matrix. With \eqref{eqsvd}, we have
\begin{equation}\label{eq4}
  x_{naive}=\sum\limits_{i=1}^{n}\frac{u_i^{T}b}{\sigma_i}v_i =
  \sum\limits_{i=1}^{n}\frac{u_i^{T}b_{true}}{\sigma_i}v_i +
  \sum\limits_{i=1}^{n}\frac{u_i^{T}e}{\sigma_i}v_i
  =x_{true}+\sum\limits_{i=1}^{n}\frac{u_i^{T}e}{\sigma_i}v_i
\end{equation}
and $\|x_{true}\|=\|A^{\dagger}b_{true}\|=
\left(\sum_{i=1}^n\frac{|u_i^Tb_{true}|^2}{\sigma_i^2}\right)^{1/2}$.

The discrete Picard condition means that, on average, the Fourier coefficient
$|u_i^{T}b_{true}|$ decays faster than $\sigma_i$, which
results in the following popular model that is used throughout Hansen's books
\cite{hansen98,hansen10} and the references therein as well as
\cite{jia18a,jia18b}:
\begin{equation}\label{picard}
  | u_i^T b_{true}|=\sigma_i^{1+\beta},\ \ \beta>0,\ i=1,2,\ldots,n,
\end{equation}
where $\beta$ is a model parameter that controls the decay rates of
$| u_i^T b_{true}|$.

The covariance matrix of the Gaussian white noise $e$
is $\eta^2 I$, the expected value $\mathcal{E}(\|e\|^2)=m \eta^2$ and  $\mathcal{E}(|u_i^Te|)=\eta,\,i=1,2,\ldots,n$, so that
$\|e\|\approx \sqrt{m}\eta$ and $|u_i^Te|\approx \eta,\ i=1,2,\ldots,n$.
\eqref{eq4} and \eqref{picard} show that, for large singular values,
$|{u_i^{T}b_{true}}|/{\sigma_i}$ is dominant relative to
$|u_i^{T}e|/{\sigma_i}$. Once
$| u_i^T b_{true}| \leq | u_i^T e|$ from some $i$ onwards, the noise
$e$ dominates $| u_i^T b|$, and the terms $\frac{| u_i^T b|}{\sigma_i}\approx
\frac{|u_i^{T}e|}{\sigma_i}$ overwhelm $x_{true}$
for small singular values and must be dampened.  Therefore, the
transition point $k_0$ is such that
\begin{equation}\label{picard1}
| u_{k_0}^T b|\approx | u_{k_0}^T b_{true}|> | u_{k_0}^T e|\approx
\eta, \ | u_{k_0+1}^T b|
\approx | u_{k_0+1}^Te|
\approx \eta;
\end{equation}
see \cite[p.42, 98]{hansen10} and \cite[p.70-1]{hansen98}.

The truncated SVD (TSVD) method \cite{hansen90,hansen98,hansen10} is a reliable
and commonly used method for solving small to modest sized
\eqref{posed}, and it solves a sequence of
problems
\begin{equation}\label{tsvd}
\min\|x\| \ \ \mbox{subject to}\ \
 \|A_kx-b\|=\min
\end{equation}
starting with $k=1$ onwards, where $A_k=U_k\Sigma_k V_k^T$
is a best rank $k$ approximation to $A$ with respect to the 2-norm
with $U_k=(u_1,\ldots,u_k)$, $V_k=(v_1,\ldots,v_k)$ and $\Sigma_k=
{\rm diag}(\sigma_1,\ldots,\sigma_k)$; it holds that
$\|A-A_k\|=\sigma_{k+1}$ \cite[p.12]{bjorck96}, and
$
x_{k}^{tsvd}=A_k^{\dagger}b
$
solves \eqref{tsvd}, called the TSVD regularized solution.
For the Gaussian white noise $e$ it is known from \cite[p.70-1]{hansen98}
and \cite[p.71,86-8,95]{hansen10} that $x_{k_0}^{tsvd}$ is
the 2-norm filtering best TSVD regularized solution of \eqref{eq1}, i.e.,
$x_{k_0}^{tsvd}$ has the minimal 2-norm error
$\|x_{true}-x_{k_0}^{tsvd}\|=\min_{k=1,2,\ldots,n}\|x_{true}-x_k^{tsvd}\|$.
The index $k$ plays the role of the regularization parameter in the TSVD method.
It has been observed and justified that $x_{k_0}^{tsvd}$ is essentially
a 2-norm filtering best possible solution of \eqref{eq1};
see \cite{hansen90b}, \cite[p.109-11]{hansen98},
\cite[Sections 4.2 and 4.4]{hansen10} and \cite{varah79}.
We refer to \cite{jia18a} for general elaborations. As a result, we can take
$x_{k_0}^{tsvd}$ as the standard reference when assessing the regularization
ability of a 2-norm filtering regularization method.

For $A$ large, the TSVD method is generally prohibitively expensive,
and only iterative regularization methods are appealing.
Krylov iterative solvers have formed a major class of methods
\cite{aster,engl00,gilyazov,hanke95,hansen98,hansen10,kirsch}.
Specifically, the CGLS method \cite{golub89,hestenes} and its mathematically
equivalent LSQR method \cite{paige82}, the CGME method
\cite{bjorck96,bjorck15,craig,hanke95,hanke01}
and the LSMR method \cite{bjorck15,chung15,fong}
have been commonly used. These methods are deterministic 2-norm filtering
regularization methods, have general regularizing
effects, and exhibit semi-convergence \cite[p.89]{natterer};
see also \cite[p.314]{bjorck96}, \cite[p.733]{bjorck15},
\cite[p.135]{hansen98} and \cite[p.110]{hansen10}: The iterates
first converge to $x_{true}$, then the
noise $e$ starts to deteriorate the iterates so that they start to diverge
from $x_{true}$ and instead converge to $x_{naive}$.
The iteration number plays the role of the regularization parameter in
iterative regularization methods.

The behavior of ill-posed problems and solvers
depends on the decay rate of $\sigma_j$. Hoffmann \cite{hofmann86}
has characterized the degree of ill-posedness
of \eqref{eq1} as follows:
If $\sigma_j=\mathcal{O}(\rho^{-j})$ with $\rho>1$,
$j=1,2,\ldots,n$, then \eqref{eq1} is severely ill-posed;
if $\sigma_j=\mathcal{O}(j^{-\alpha})$, then \eqref{eq1}
is mildly or moderately ill-posed for $\frac{1}{2}<\alpha\le1$ or $\alpha>1$.
This definition has been widely used \cite{aster,engl00,hansen98,hansen10}.
The requirement $\alpha>\frac{1}{2}$ does not
appear in \cite{hofmann86} and is explicitly added in
\cite{huangjia,jia18a}, which is always met
for a linear compact operator equation \cite{hanke93,hansen98}.

Hanke and Hansen \cite{hanke93} address that
a strict proof of the regularizing properties of conjugate gradients is
extremely difficult; see also \cite{hansen07}. The regularizing effects
of CGLS, LSQR and CGME have been intensively studied; see, e.g.,
and have been intensively studied
\cite{aster,eicke,firro97,gilyazov,hanke95,hanke01,hansen98,hansen10,hps16,
hps09,huangjia,jia18a,jia18b,paige06,vorst90}.
It has long been known (cf. \cite{hanke93,hansen98,hansen07,hansen10})
that if the singular values of the projection matrices involved
in LSQR, called the Ritz values, approximate
the large singular values in natural order then LSQR has
the same regularization ability as the TSVD method, that is, the two
methods can compute 2-norm filtering best regularized solutions with
the same accuracy. As we will see clearly, the same results hold for CGME and
LSMR when the singular values of projection matrices approximate the large
singular values of $A$ and $A^TA$ in this order, respectively.

If a 2-norm filtering regularized solution of \eqref{eq1} is
as accurate as $x_{k_0}^{tsvd}$, it is called a 2-norm filtering
best possible regularized solution. If the 2-norm filtering regularized
solution by a regularization method at semi-convergence is
such a best possible one, then the solver is said to have
the {\em full} regularization. Otherwise, the solver has
only the {\em partial} regularization. This definition is introduced in
\cite{huangjia,jia18a}. In terms of it, a fundamental question posed
in \cite{huangjia,jia18a} is:
{\em  Do CGLS, LSQR, CGME and LSMR have the full or partial regularization for
severely, moderately and mildly ill-posed problems?} Actually, this question
has been receiving high attention for CGLS and LSQR.

For the cases that $\sigma_i$ are simple,
the author in \cite{jia18a} has given accurate estimates for the 2-norm
distances between the underlying $k$ dimensional Krylov subspace and the
$k$ dimensional dominant right singular subspace $span\{V_k\}$
of $A$ for severely, moderately and mildly ill-posed problems.
On the basis of \cite{jia18a}, the author in \cite{jia18b} has proved
that, for LSQR, the $k$ Ritz values converge to the $k$
large singular values of $A$ in natural order and Lanczos bidiagonalization
always generates a near best rank $k$ approximation until $k=k_0$
for severely and moderately
ill-posed problems with suitable $\rho>1$ and $\alpha>1$, meaning
that LSQR and CGLS have the full regularization. However, if
such desired properties fail to hold, it has been theoretically
unknown if LSQR has the full or partial regularization. Nevertheless,
numerical experiments on many
ill-posed problems have demonstrated that LSQR always
has the full regularization \cite{jia18a,jia18b}.

In this paper, we analyze the regularization of CGME and LSMR under the
assumption that all the singular values $\sigma_i$ are simple.
We establish a number of results, and prove that the regularization
ability of CGME is generally inferior to that of LSQR,
that is, the 2-norm filtering best regularized solutions obtained
by CGME at semi-convergence
are generally less accurate than those obtained by LSQR.
Specifically, we derive the filtered SVD expansion of CGME iterates,
by which we prove that the semi-convergence of CGME always occurs no later
than that of LSQR and can be much earlier than the latter.
In the meantime, we show how to extract a rank $k$
approximation from the rank $k+1$ approximation to $A$ generated
in CGME at iteration $k$, which is
as accurate as the rank $k$ approximation in LSQR. Exploiting such
rank $k$ approximation,
we propose a modified CGME (MCGME) method whose regularization ability is shown
to be very comparable to that of LSQR. For LSMR, we present a number of results
and prove that its regularization ability is as good as that of LSQR
and the two methods compute the 2-norm filtering best regularized solutions
with essentially the same accuracy. We also show that the
semi-convergence of LSMR always occurs no sooner than that of LSQR.

As a windfall, making of our analysis approach used for CGME,
we improve a fundamental bound, Theorem 9.3
presented in Halko {\em et al.} \cite{halko11}, for the
accuracy of the truncated rank $k$ approximation to $A$
generated by randomized algorithms, which have formed a highly intensive
topic and have been used in numerous disciplines
over the years. As remarked by Halko {\em et al.}
in \cite{halko11} (cf. Remark 9.1 there), their bound
appears ``{\em conservative, but a complete theoretical
understanding lacks.}'' Our new bounds for the approximation
accuracy are not only unconditionally sharper than theirs but also
can reveal how the truncation step damages the accuracy
of the rank $k$ approximation.

The paper is organized as follows. In Section~\ref{methods}, we
review LSQR, CGME and LSMR. In Section~\ref{lsqr}, we briefly
state some results on LSQR in \cite{jia18a,jia18b} and
take LSQR as reference to assess the regularization ability of CGME
and LSMR. In Section~\ref{cgme}, we derive a number of regularization
properties of CGME and propose the MCGME method.
In Section~\ref{randomappro}, we consider the accuracy
of the truncated rank $k$ randomized approximation \cite{halko11} and present sharper
bounds. In Section~\ref{lsmr}, we study the regularization ability of LSMR.
In Section~\ref{numer}, we report
numerical experiments to confirm our theory. We conclude
the paper in Section~\ref{conclu}.

Throughout the paper, we denote by
$\mathcal{K}_{k}(C, w)= span\{w,Cw,\ldots,C^{k-1}w\}$
the $k$ dimensional Krylov subspace generated
by the matrix $\mathit{C}$ and the vector $\mathit{w}$, and by
the bold letter $\mathbf{0}$ the zero matrix with orders clear from the context.

\section{The LSQR, CGME and LSMR algorithms}\label{methods}

These three algorithms are all based on the Lanczos bidiagonalization process,
which computes two orthonormal bases $\{q_1,q_2,\dots,q_k\}$ and
$\{p_1,p_2,\dots,p_{k+1}\}$  of $\mathcal{K}_{k}(A^{T}A,A^{T}b)$ and
$\mathcal{K}_{k+1}(A A^{T},b)$  for $k=1,2,\ldots,n$,
respectively. We describe the process as Algorithm 1.

{\bf Algorithm 1: \  $k$-step Lanczos bidiagonalization process}
\begin{remunerate}
\item Take $ p_1=b/\|b\| \in \mathbb{R}^{m}$, and define $\beta_1{q_0}=\mathbf{0}$.

\item For $j=1,2,\ldots,k$
 \begin{romannum}
  \item
  $r = A^{T}p_j - \beta_j{q_{j-1}}$
  \item $\alpha_j = \|r\|;q_j = r/\alpha_j$
  \item
   $   z = Aq_j - \alpha_j{p_{j}}$
  \item
  $\beta_{j+1} = \|z\|;p_{j+1} = z/\beta_{j+1}.$
   \end{romannum}
\end{remunerate}

Algorithm 1 can be written in the matrix form
\begin{align}
  AQ_k&=P_{k+1}B_k,\label{eqmform1}\\
  A^{T}P_{k+1}&=Q_{k}B_k^T+\alpha_{k+1}q_{k+1}(e_{k+1}^{(k+1)})^{T},\label{eqmform2}
\end{align}
where $e_{k+1}^{(k+1)}$ denotes the $(k+1)$-th canonical basis vector of
$\mathbb{R}^{k+1}$, $P_{k+1}=(p_1,p_2,\ldots,p_{k+1})$,
$Q_k=(q_1,q_2,\ldots,q_k)$ and
\begin{equation}\label{bk}
  B_k = \left(\begin{array}{cccc} \alpha_1 & & &\\ \beta_2 & \alpha_2 & &\\ &
  \beta_3 &\ddots & \\& & \ddots & \alpha_{k} \\ & & & \beta_{k+1}
  \end{array}\right)\in \mathbb{R}^{(k+1)\times k}.
\end{equation}
It is known from \eqref{eqmform1} that
\begin{equation}\label{Bk}
B_k=P_{k+1}^TAQ_k.
\end{equation}
Algorithm 1 cannot break down before step $n$ when
$\sigma_i,\ i=1,2,\ldots,n$, are simple since $b$ is supposed to have
nonzero components in the directions of $u_i,\ i=1,2,\ldots,n$.
The singular values $\theta_i^{(k)},\ i=1,2,\ldots,k$ of $B_k$,
called the Ritz values of $A$ with
respect to the left and right subspaces $span\{P_{k+1}\}$ and $span\{Q_k\}$,
are all simple.

Write
$\mathcal{V}_k^R=\mathcal{K}_k(A^TA,A^Tb)$ and $\beta_1=\|b\|$.
At iteration $k$, LSQR \cite{paige82} solves
$$
\|Ax_k^{lsqr}-b\|=\min_{x\in \mathcal{V}_k^R}
\|Ax-b\|
$$
for the iterate
\begin{equation}\label{yk}
x_k^{lsqr}=Q_ky_k^{lsqr} \ \ \mbox{with}\ \
y_k^{lsqr}=\arg\min\limits_{y\in \mathbb{R}^{k}}\|B_ky-\beta_1 e_1^{(k+1)}\|
  =\beta_1 B_k^{\dagger} e_1^{(k+1)},
\end{equation}
where $e_1^{(k+1)}$ is the first canonical basis vector of $\mathbb{R}^{k+1}$,
and $\|Ax_k^{lsqr}-b\|=\|B_ky_k^{lsqr}-\beta_1 e_1^{(k+1)}\|$
decreases monotonically with respect to $k$.

CGME \cite{bjorck15,hanke95,hanke01,hps16,hps09} is the CG method
implicitly applied to $\min\|AA^Ty-b\|$ or $AA^Ty=b$ with $x=A^Ty$, and it
solves the problem
$$
\|x_{naive}-x_k^{cgme}\|=\min_{x\in \mathcal{V}_k^R}\|x_{naive}-x\|
$$
for the iterate $x_k^{cgme}$. The error norm
$\|x_{naive}-x_k^{cgme}\|$ decreases monotonically with respect to $k$.
Let $\bar{B}_k\in \mathbb{R}^{k\times k}$ be the
matrix consisting of the first $k$ rows of $B_k$, i.e.,
\begin{equation}\label{bbar}
\bar{B}_k=P_k^TAQ_k.
\end{equation}
Then the CGME iterate
\begin{equation}\label{ykcgme}
x_k^{cgme}=Q_ky_k^{cgme} \ \ \mbox{with} \ \
y_k^{cgme}=\beta_1 \bar{B}_k^{-1}e_1^{(k)}
\end{equation}
and $\|Ax_k^{cgme}-b\|=\beta_{k+1}|(e_k^{(k)})^T y_k^{cgme}|$
with $e_k^{(k)}$ the $k$-th canonical vector of $\mathbb{R}^{k+1}$.

LSMR \cite{bjorck15,fong} is mathematically equivalent to MINRES \cite{paige75}
applied to the normal equation $A^TAx=A^Tb$ of \eqref{eq1}, and it solves
$$
\|A^T(b-A x_k^{lsmr})\|=\min_{x\in \mathcal{V}_k^R}\|A^T(b-A x)\|
$$
for the iterate $x_k^{lsmr}$.
The residual norm $\|A^T(b-Ax_k^{lsmr})\|$ of the normal equation decreases
monotonically with respect to $k$, and the iterate
\begin{equation}\label{yklsmr}
x_k^{lsmr}=Q_ky_k^{lsmr} \ \ \mbox{with}\ \
y_k^{lsmr}=\arg\min\limits_{y\in \mathbb{R}^{k}}\|(B_k^TB_k,\alpha_{k+1}
\beta_{k+1}e_k^{(k)})^Ty-\alpha_1\beta_1 e_1^{(k+1)}\|.
\end{equation}

\section{Some results on LSQR in \cite{jia18a,jia18b}} \label{lsqr}

From  $\beta_1 e_1^{(k+1)}=P_{k+1}^T b$ and \eqref{yk} we have
\begin{equation}\label{xk}
x_k^{lsqr}=Q_k B_k^{\dagger} P_{k+1}^Tb,
\end{equation}
which is the minimum 2-norm
solution to the problem that perturbs $A$ in \eqref{eq1} to
its rank $k$ approximation $P_{k+1}B_k Q_k^T$. Recall that $\|A-A_k\|=\sigma_{k+1}$.
Analogous to \eqref{tsvd}, LSQR now solves a sequence of problems
\begin{equation}\label{lsqrreg}
\min\|x\| \ \ \mbox{ subject to }\ \ \|P_{k+1}B_kQ_k^Tx-b\|=\min
\end{equation}
for $x_k^{lsqr}$ starting with $k=1$
onwards, where $A$ in \eqref{eq1} is replaced by a rank $k$ approximation
$P_{k+1}B_kQ_k^T$ of it. Therefore, if $P_{k+1}B_k Q_k^T$ is a near best rank $k$
approximation to $A$ with an approximate accuracy $\sigma_{k+1}$ and
the singular values $\theta_i^{(k)},\
i=1,2,\ldots,k$ of $B_k$ approximate the $k$
large $\sigma_i$ in natural order for $k=1,2,\ldots,k_0$, then
LSQR has the same regularization ability as the TSVD method
and thus has the full regularization. See \cite{jia18a}
for more elaborations.

The analysis on the TSVD method and the Tikhonov regularization
method \cite{hansen98,hansen10} shows that the core requirement on
a regularization method is to
acquire the $k_0$ dominant SVD components of $A$ and meanwhile
suppress the remaining $n-k_0$ SVD components. Therefore, the more accurate
the rank $k$ approximation is to $A$ and the better approximations are
the $k$ non-zero singular values of a projection matrix
to some of the $k_0$ large singular values of $A$, the better
regularization ability of the method has, so that
the best regularized solution obtained by it is more accurate.

Define
\begin{equation}\label{gammak}
\gamma_k^{lsqr} = \|A-P_{k+1}B_kQ_k^T\|,
\end{equation}
which measures the accuracy of the rank $k$ approximation $P_{k+1}B_kQ_k^T$
to $A$ involved in LSQR. Since
the best rank $k$ approximation $A_k$ satisfies $\|A-A_k\|=\sigma_{k+1}$,
we have
$$
\gamma_k^{lsqr}\geq \sigma_{k+1}.
$$

The author in \cite{jia18b} introduces the definition of
a near best rank $k$ approximation to $A$: For LSQR,
$P_{k+1}B_kQ_k^T$ is called a near best rank $k$ approximation to $A$ if
$\gamma_k^{lsqr}$ is closer to $\sigma_{k+1}$ than to $\sigma_k$:
\begin{equation}\label{near}
\sigma_{k+1}\leq \gamma_k^{lsqr}<\frac{\sigma_k+\sigma_{k+1}}{2}.
\end{equation}

Based on the accurate estimates established in \cite{jia18a}
for the 2-norm distances between the underlying Krylov subspace
$\mathcal{V}_k^R$ and the $k$ dimensional dominant right singular
subspace $span\{V_k\}$ for severely, moderately and
mildly ill-posed problems,
the author \cite{jia18b} has derived accurate estimates for $\gamma_k^{lsqr}$
and a number of approximation properties of $\theta_i^{(k)},\ i=1,2,\ldots,k$
for the three kinds of ill-posed problems.
The results have shown that, for severely and moderately ill-posed problems
with  for suitable $\rho>1$ and $\alpha>1$ and for $k=1,2,\ldots,k_0$,
$P_{k+1}B_kQ_k^T$ must be a near best rank $k$ approximation to $A$, and
the $k$ Ritz values $\theta_i^{(k)}$ approximate
the large singular values $\sigma_i$ of $A$
in natural order. This means
that LSQR has the full regularization for these two kinds of problems with
suitable $\rho>1$ and $\alpha>1$. However, for moderately ill-posed problems
with $\alpha>1$ not enough and mildly ill-posed problems,
$P_{k+1}B_kQ_k^T$ is generally not a near best rank $k$ approximation,
and the $k$ Ritz values $\theta_i^{(k)}$ do not approximate
the large singular values of $A$ in natural order for some $k\leq k^*$.

In particular, the author \cite[Theorem 5.1]{jia18b} has proved
the following three results:
\begin{align}
\gamma_k^{lsqr}&=\|G_k\| \label{gk}
\end{align}
with
\begin{align}\label{gk1}
G_k&=\left(\begin{array}{cccc}
\alpha_{k+1} & & & \\
\beta_{k+2}& \alpha_{k+2} & &\\
&
  \beta_{k+3} &\ddots & \\& & \ddots & \alpha_{n} \\
  & & & \beta_{n+1}
  \end{array}\right)\in \mathbb{R}^{(n-k+1)\times (n-k)},
\end{align}
\begin{align}
\alpha_{k+1}&<\gamma_k^{lsqr},\
\beta_{k+2}<\gamma_k^{lsqr},\ k=1,2,\ldots,n-1,\label{alphagamma}\\
\gamma_{k+1}^{lsqr}&<\gamma_k^{lsqr},\ k=1,2,\ldots,n-2.\label{monto}
\end{align}
These notation and results will be used later.

\section{The regularization of CGME}\label{cgme}

Note that $P_k^Tb=\beta_1e_1^{(k)}$. We obtain
\begin{equation}\label{cgmesolution}
x_k^{cgme}=Q_k\bar{B}_k^{-1}P_k^Tb.
\end{equation}
Therefore, analogous to \eqref{tsvd} and \eqref{lsqrreg}, CGME solves
a sequence of problems
\begin{equation}\label{cgmereg}
\min\|x\| \ \ \mbox{ subject to }\ \ \|P_k\bar{B}_kQ_k^T x-b\|=\min
\end{equation}
for the regularized solution $x_k^{cgme}$ starting
with $k=1$ onwards, where $A$ in \eqref{eq1} is replaced by
a rank $k$ approximation $P_k\bar{B}_kQ_k^T$ of it.

Just as LSQR, if $P_k\bar{B}_kQ_k^T$ is a near best rank $k$ approximation to
$A$ and the $k$ singular values of $\bar{B}_k$ approximate the
large ones of $A$ in natural order for $k=1,2,\ldots,k_0$, then
CGME has the full regularization.

By \eqref{eqmform1}, \eqref{eqmform2} and \eqref{Bk},
the rank $k$ approximation involved in LSQR is
\begin{equation}\label{lsqrappr}
P_{k+1}B_kQ_k^T=AQ_kQ_k^T.
\end{equation}
By \eqref{gammak}, we have
$
\gamma_k^{lsqr}=\|A(I-Q_kQ_k^T)\|.
$
For CGME, by \eqref{eqmform2} and \eqref{bbar}, we obtain
\begin{align}
P_{k+1}P_{k+1}^TA &= P_{k+1}(B_kQ_k^T+\alpha_{k+1}e_{k+1}^{(k+1)}q_{k+1}^T)\notag\\
&=P_{k+1}(B_k, \alpha_{k+1}e_{k+1}^{(k+1)})Q_{k+1}^T \notag\\
&=P_{k+1}\bar{B}_{k+1}Q_{k+1}^T. \label{leftbidiag}
\end{align}
Therefore, $x_k^{cgme}$ is the solution to \eqref{cgmereg}
in which the rank $k$ approximation to $A$ is
$P_k\bar{B}_kQ_k^T=P_kP_k^TA$, whose approximation accuracy
is
\begin{equation}\label{cgmeacc}
\gamma_k^{cgme}=\|A-P_k\bar{B}_kQ_k^T\|=\|(I-P_kP_k^T)A\|.
\end{equation}

\begin{theorem}\label{cgmeappr}
For the rank $k$ approximations $P_kP_k^TA=P_k\bar{B}_kQ_k^T$ to $A$,
$k=1,2,\ldots,n-1$, with the definition $\gamma_0^{lsqr}=\|A\|$ we have
\begin{eqnarray}
&\gamma_k^{lsqr}<\gamma_k^{cgme}< \gamma_{k-1}^{lsqr}&, \label{cgmelowup} \\
&\gamma_{k+1}^{cgme}< \gamma_k^{cgme}.&  \label{mono}
\end{eqnarray}
\end{theorem}

{\em Proof. }
We give two proofs of the upper bound in \eqref{cgmelowup}. The first is as follows.
Since $P_{k+1}P_{k+1}^T(I-P_{k+1}P_{k+1}^T)=\mathbf{0}$,
from \eqref{eqmform2} we obtain
\begin{align}
  (\gamma_k^{lsqr})^2 &=\|A-P_{k+1}B_kQ_k^T\|^2 \notag\\
  &=\|P_{k+1}P_{k+1}^TA-P_{k+1}B_kQ_k^T+(I-P_{k+1}P_{k+1}^T)A\|^2\notag\\
  &=\max_{\|y\|=1}\|\left(
  (P_{k+1}P_{k+1}^TA-P_{k+1}B_kQ_k^T)+(I-P_{k+1}P_{k+1}^T)A\right)y\|^2\notag\\
   &=\max_{\|y\|=1} \|P_{k+1}P_{k+1}^T(P_{k+1}P_{k+1}^TA-
           P_{k+1}B_kQ_k^T)y+(I-P_{k+1}P_{k+1}^T)Ay\|^2 \notag\\
           &=\max_{\|y\|=1}\left(\|P_{k+1}P_{k+1}^T(P_{k+1}P_{k+1}^TA-
           P_{k+1}B_kQ_k^T)y\|^2+\|(I-P_{k+1}P_{k+1}^T)Ay\|^2\right)\notag\\
            &=\max_{\|y\|=1}\left(\|P_{k+1}(P_{k+1}^TA-B_kQ_k^T)y\|^2+
           \|(I-P_{k+1}P_{k+1}^T)Ay\|^2\right)\notag\\
           &=\max_{\|y\|=1}\left(\|(P_{k+1}^TA-B_kQ_k^T)y\|^2+
           \|(I-P_{k+1}P_{k+1}^T)Ay\|^2\right)\notag\\
            &=\max_{\|y\|=1}\left(\alpha_{k+1}^2|(e_{k+1}^{(k+1)})^Ty|^2+
           \|(I-P_{k+1}P_{k+1}^T)Ay\|^2\right)\notag\\
            &>
            \max_{\|y\|=1} \|(I-P_{k+1}P_{k+1}^T)Ay\|^2 \notag\\
            &=\|(I-P_{k+1}P_{k+1}^T)A\|^2=(\gamma_{k+1}^{cgme})^2,\notag
\end{align}
which is the upper bound in \eqref{cgmelowup} by replacing the index $k+1$ with $k$.

Taking $k=n$ in \eqref{Bk} and augmenting $P_{n+1}$ such that
$P=(P_{n+1},\widehat{P})\in \mathbb{R}^{m\times m}$ is orthogonal, we have
\begin{equation}\label{fulllb}
P^TAQ_n=\left(\begin{array}{c}
B_n\\
\mathbf{0}
\end{array}
\right),
\end{equation}
where all the entries $\alpha_i$ and $\beta_{i+1}$, $i=1,2,\ldots,n$, of
$B_n$ are positive, and $Q_n\in \mathbb{R}^{n\times n}$ is orthogonal.
Then by the orthogonal invariance of the 2-norm we obtain
\begin{equation}\label{cgmeapp}
\gamma_k^{cgme}=\|A-P_k\bar{B}_kQ_k^T\|=\|P^T(A-P_k\bar{B}_kQ_k^T)Q_n\|
=\|(\beta_{k+1}e_1,G_k)\|
\end{equation}
with $G_k$ defined by \eqref{gk1}. It is straightforward to justify
that the singular values of $G_k\in \mathbb{R}^{(n-k+1)\times (n-k)}$ {\em strictly}
interlace those of $(\beta_k e_1,G_k)\in \mathbb{R}^{(n-k+1)\times (n-k+1)}$
by noting that $(\beta_{k+1}e_1,G_k)^T (\beta_{k+1}e_1,G_k)$ is
an {\em unreduced} symmetric tridiagonal matrix, from which and
$\|G_k\|=\gamma_k^{lsqr}$ the lower bound of \eqref{cgmelowup} follows.

Based on \eqref{cgmeapp}, we can also give the second proof of the upper bound
in \eqref{cgmelowup}.
Observe from \eqref{gk1} that $(\beta_{k+1}e_1,G_k)$ is the matrix deleting
the first row of $G_{k-1}$. Applying the strict interlacing property of
singular values to $(\beta_{k+1}e_1,G_k)$ and $G_{k-1}$, we obtain
$\gamma_{k-1}^{lsqr}=\|G_{k-1}\|>\|(\beta_{k+1}e_1,G_k)\|=\gamma_k^{cgme}$,
which yields the upper bound of \eqref{cgmelowup}.

From \eqref{cgmeapp}, notice that $(\beta_{k+2}e_1,G_{k+1})$ is the matrix
deleting the first row of $(\beta_{k+1}e_1,G_k)$ and the first
column, which is {\em zero}, of the resulting matrix.
Applying the strict interlacing property of singular values
to $(\beta_{k+2}e_1,G_{k+1})$ and $(\beta_{k+1}e_1,G_k)$ establishes
\eqref{mono}.
\qquad\endproof

\eqref{cgmelowup}
indicates that $P_kP_k^TA=P_k\bar{B}_kQ_k^T$ is definitely a less accurate
rank $k$ approximation to $A$ than $AQ_kQ_k^T=P_{k+1}B_kQ_k^T$ in LSQR.
\eqref{mono} shows the strict monotonic decreasing property of $\gamma_k^{cgme}$.
Moreover, keep in mind that $\gamma_k^{lsqr}\geq \sigma_{k+1}$. Then
a combination of it and the results in Section~\ref{lsqr}
indicates that, unlike $P_{k+1}B_kQ_k^T$ in LSQR,
there is no guarantee that $P_k\bar{B}_kQ_k^T$ is a near best
rank $k$ approximation to $A$ even for severely and moderately
ill-posed problems, because $\gamma_k^{cgme}$ simply lies
between $\gamma_k^{lsqr}$ and $\gamma_{k-1}^{lsqr}$ and
there do not exist any sufficient conditions on
$\rho>1$ and $\alpha>1$ that enforce $\gamma_k^{cgme}$ to be closer
to $\gamma_k^{lsqr}$, let alone closer to $\sigma_{k+1}$.
Therefore, based on the accuracy of the rank $k$ approximations
in CGME and LSQR, we come to the conclusion that
the regularization ability of CGME cannot be
superior and is generally inferior to that of LSQR. Furthermore,
since there is no guarantee that $P_k\bar{B}_kQ_k^T$ is a near best
rank $k$ approximation for severely and moderately ill-posed
problems with suitable $\rho>1$ and $\alpha>1$, CGME
may not have the full regularization for these two kinds of problems.

In the following we investigate the approximation behavior of the $k$ singular
values $\bar{\theta}_i^{(k)}$ of $\bar{B}_k,\ k=1,2,\ldots,n$.
Before proceeding, it is necessary to have a closer look
at Algorithm 1 and distinguish some subtleties when $A$ is rectangular,
i.e., $m>n$, and square, i.e., $m=n$, respectively.

Keep in mind that Algorithm 1 does not break down
before step $n$. For the rectangular case $m>n$, Algorithm 1 is exactly what
is presented
there, all the $\alpha_k$ and $\beta_{k+1}$ are positive, $k=1,2,\ldots,n$,
and we generate $P_{n+1}$ and $Q_n$ at step $n$ and $\alpha_{n+1}=\beta_{n+2}=0$.
As a consequence, by definition \eqref{leftbidiag}, we have
\begin{equation}\label{bn1}
\bar{B}_{n+1}=(B_n,\alpha_{n+1}e_{n+1}^{(n+1)})=(B_n,\mathbf{0}).
\end{equation}
It is known from \eqref{fulllb} that the singular values of $B_n$ are
identical to the singular values
$\sigma_i,\ i=1,2,\ldots,n$ of $A$. Therefore, the $n+1$
singular values of $\bar{B}_{n+1}$ are $\sigma_i,\,i=1,2,\ldots,n$
and zero.

For the square case $m=n$,
however, we must have $\beta_{n+1}=0$, that is, the last row of $B_n$
is zero; otherwise, we would obtain an $n\times (n+1)$ orthonormal matrix $P_{n+1}$,
which is impossible since $P_n$ is already an orthogonal matrix.
After Algorithm 1 is run to completion, we have
$$
\bar{B}_n=P_n^TAQ_n,
$$
whose singular values $\bar{\theta}_i^{(n)}=\sigma_i,\,i=1,2,\ldots,n$.

By the definition \eqref{leftbidiag} of $\bar{B}_k$,
from \eqref{eqmform2} and the above description,
for both the rectangular and square cases we obtain
\begin{equation}\label{aat}
P_k^TAA^TP_k=\bar{B}_k \bar{B}_k^T,\,k=1,2,\ldots,n^*,
\end{equation}
with $n^*=n+1$ for $m>n$ and $n^*=n$ for $m=n$,
which are unreduced symmetric tridiagonal matrices.
For $m=n$, the eigenvalues of $AA^T$ are just
$\sigma_i^2,\ i=1,2,\ldots,n$, all of which are simple and positive;
for $m>n$, the eigenvalues of $AA^T$ are $\sigma_i^2,\,
i=1,2,\ldots,n$ plus $m-n$ zeros, denoted by $\sigma_{n+1}^2=\cdots=\sigma_m^2=0$
for our later use. Therefore, by the definition of $n^*$,
the eigenvalues of $\bar{B}_{n^*}\bar{B}_{n^*}^T$ are
$\sigma_i^2,\ i=1,2,\ldots,n^*$.

Notice that
$\bar{B}_k \bar{B}_k^T$ is nothing but the projection matrix of
$AA^T$ onto the $k$ dimensional Krylov subspace $\mathcal{K}_k(AA^T,b)$.
More precisely, $\bar{B}_k \bar{B}_k^T$  is generated by
the $k$-step symmetric Lanczos tridiagonalization process
applied to $AA^T$ starting with $p_1=b/\|b\|$, and the eigenvalues
of $\bar{B}_k \bar{B}_k^T$ generally approximate extreme eigenvalues of $AA^T$;
see, e.g., \cite{bjorck96,bjorck15,parlett} for details.
Particularly, the smallest eigenvalue $(\bar{\theta}_k^{(k)})^2$
of $\bar{B}_k \bar{B}_k^T$ generally converges to the smallest
eigenvalue $\sigma_{n^*}^2$ of $AA^T$ as $k$ increases, which is
$\sigma_{n+1}^2=0$ for $m>n$ and $\sigma_n^2>0$ for $m=n$.
In contrast, for $B_k$, its smallest singular value
$\theta_k^{(k)}> \sigma_n$ unconditionally until $\theta_n^{(n)}=\sigma_n$.

We next give a number of close relationships between
$\bar{\theta}_i^{(k)}$ and $\theta_i^{(k)}$ as well as between them and
the singular values $\sigma_i$ of $A$, which are crucial to compare the
regularizing effects of CGME with those of LSQR.

\begin{theorem}\label{interlace}
Denote by $\bar{\theta}_i^{(k)}$ and $\theta_i^{(k)},\ i=1,2,\ldots,k$
the singular values of $\bar{B}_k$ and $B_k$, respectively,
labeled in decreasing order. Then
\begin{align}
\theta_1^{(k)}&>\bar{\theta}_1^{(k)}>\theta_2^{(k)}>
\bar{\theta}_2^{(k)}> \cdots >\theta_k^{(k)}> \bar{\theta}_k^{(k)}, \
k=1,2,\ldots,n-1. \label{secondinter}
\end{align}
Moreover,
\begin{align}
\sigma_n&<\bar{\theta}_k^{(k)}<\theta_k^{(k)}<\sigma_k,\ k=1,2,\ldots,n-1
\label{thetak}
\end{align}
for $m=n$ and
\begin{align}
\sigma_n&<\theta_k^{(k)}<\sigma_k,\ k=1,2,\ldots,n-1,\label{thetak1} \\
0&<\bar{\theta}_k^{(k)}<\theta_k^{(k)}<\sigma_k,\ k=1,2,\ldots,n-1 \label{thetak2}
\end{align}
for $m>n$.
\end{theorem}

{\em Proof.}
Observe that $\bar{B}_k$ consists of the first $k$ rows of $B_k$ and
all the $\alpha_k$ and $\beta_{k+1}$ are positive for $k=1,2,\ldots,n-1$.
Applying the strict interlacing property of singular values to $\bar{B}_k$
and $B_k$, we obtain \eqref{secondinter}.

Note that, for $A$ both rectangular and square,
we have $\theta_i^{(n)}=\sigma_i,\ i=1,2,\ldots,n$. Since $B_k$ consists
of the first $k$ columns of $B_n$ and deletes the last $n-k$ zero rows of
the resulting matrix, applying the strict interlacing property of singular
values to $B_k$ and $B_n$ (cf. \cite[p.198, Corollary 4.4]{stewartsun}),
for $k=1,2,\ldots,n-1$ we have
\begin{equation}\label{interbar}
\sigma_{n-k+i}<\theta_i^{(k)}<\sigma_i,\ i=1,2,\ldots,k.
\end{equation}

Observe that
$\bar{B}_k\bar{B}_k^T,\, k=1,2,\ldots,n-1,$ are the $k\times k$ leading
principal matrices of $\bar{B}_{n^*}\bar{B}_{n^*}^T$, whose eigenvalues are
$\sigma_i^2,\ i=1,2,\ldots,n^*$, and they are unreduced
symmetric tridiagonal matrices. Applying the strict interlacing
property of eigenvalues to $\bar{B}_k \bar{B}_k^T$ and $\bar{B}_{n^*}\bar{B}_{n^*}^T$,
for $k=1,2,\ldots,n-1$ we obtain
$$
\sigma_{n^*-k+i}^2<(\bar{\theta}_i^{(k)})^2<\sigma_i^2,\ i=1,2,\ldots,k,
$$
from which and the definition of $n^*$ it follows that
$$
\sigma_n<\bar{\theta}_k^{(k)}<\sigma_k
$$
for $m=n$ and
$$
0=\sigma_{n+1}<\bar{\theta}_k^{(k)}<\sigma_k
$$
for $m>n$.
The above, together with \eqref{interbar} and \eqref{secondinter}, yields
\eqref{thetak}--\eqref{thetak2}.
\qquad\endproof

From Section~\ref{lsqr}, \eqref{thetak} and \eqref{thetak2} indicate that,
unlike the $k$ singular values $\theta_i^{(k)}$ of $B_k$, which have been
proved to interlace the first $k+1$ large ones of $A$ and approximate
the first $k$ ones in natural order
for the severely or moderately ill-posed problems for suitable
$\rho>1$ or $\alpha>1$ \cite{jia18b}, the
lower bound for $\bar{\theta}_k^{(k)}$ is simply $\sigma_n$ for $m=n$
and zero for $m>n$, respectively,
and there does not exist a better lower bound for it. This implies that
$\bar{\theta}_k^{(k)}$ may be much smaller than $\sigma_{k+1}$ and it
can be as small as $\sigma_n$ for $m=n$ and arbitrarily small for $m>n$,
independent of $\rho$ or $\alpha$.
In other words, the size of $\rho$ or $\alpha$ has no
intrinsic effects on the size of $\bar{\theta}_k^{(k)}$, and
cannot make $\bar{\theta}_k^{(k)}$ lie between $\sigma_{k+1}$
and $\sigma_k$ by choosing $\rho$ or $\alpha$,
that is, the regularizing effects of CGME have intrinsic
indeterminacy for severely and moderately ill-posed problems, independent
of the size of $\rho$ and $\alpha$. Therefore, CGME
may or may not have the full regularization for these
two kinds of problems. On the other hand,
even if the $\bar{\theta}_i^{(k)}$ approximate the first $k$ large singular
values $\sigma_i$ in natural order, they are less accurate than
the $k$ singular values $\theta_i^{(k)}$
of $B_k$ because of \eqref{thetak} and \eqref{thetak2}. Consequently, since
the $\theta_i^{(k)}$ are always correspondingly larger than
the $\bar{\theta}_i^{(k)}$, the regularization ability of CGME
cannot be superior and is generally inferior to that of the LSQR.


A final note is that, unlike for $m=n$,
CGME may be at risk for $m>n$ since the $\bar{\theta}_k^{(k)}$ converges
to {\em zero} other than $\sigma_n$ as $k$ increases and can be arbitrarily
small, which causes that the projected problem
$\bar{B}_ky_k^{cgme}=\beta_1 e_1^{(k)}$
may even be worse conditioned than \eqref{eq1} and $\|x_k^{cgme}\|
=\|y_k^{cgme}\|$ may be unbounded as $k$ increases
and bigger than $\|x_{naive}\|$ for a given \eqref{eq1}.

In what follows we establish more results on the regularization of CGME
and get more insight into it.
It is known, e.g., \cite[p.146]{hansen98} that the LSQR iterate $x_k^{lsqr}$
takes the following filtered SVD expansion:
\begin{equation}\label{eqfilter2}
  x_k^{lsqr}=\sum\limits_{i=1}^nf_i^{(k,lsqr)}\frac{u_i^{T}b}{\sigma_i}v_i,\
  k=1,2,\ldots,n,
\end{equation}
where the filters
\begin{equation}\label{filterlsqr}
f_i^{(k,lsqr)}=1-\prod\limits_{j=1}^k\frac{(\theta_j^{(k)})^2-\sigma_i^2}
{(\theta_j^{(k)})^2},\ i=1,2,\ldots,n.
\end{equation}
These results have been extensively used to study the regularizing
effects of LSQR; see, e.g., \cite{hansen98,hansen07,jia18a}.
We now prove that the CGME iterate $x_k^{cgme}$ also takes a filtered
SVD expansion similar to \eqref{eqfilter2} and \eqref{filterlsqr},
but its proof is much more involved than that
of \eqref{eqfilter2} and \eqref{filterlsqr}.

\begin{theorem}\label{theocgme}
The CGME iterate $x_k^{cgme}$ has the filtered SVD expansion
\begin{equation}\label{cgmeexpr}
  x_k^{cgme}=\sum\limits_{i=1}^nf_i^{(k,cgme)}\frac{u_i^{T}b}{\sigma_i}v_i,\
  k=1,2,\ldots,n,
\end{equation}
where the filters
\begin{equation}\label{filter}
f_i^{(k,cgme)}=1-\prod\limits_{j=1}^k\frac{(\bar{\theta}_j^{(k)})^2-\sigma_i^2}
{(\bar{\theta}_j^{(k)})^2},\ i=1,2,\ldots,n.
\end{equation}
\end{theorem}

{\em Proof.}
Let $y_{naive}=(AA^T)^{\dagger}b$ be the minimal 2-norm
solution to
$\min_{y} \|AA^Ty-b\|$. Recall Algorithm 1. For this minimization problem,
starting with $y_0^{cgme}=\mathbf{0}$, at iteration $k$ the CG method
extracts $y_k^{cgme}$ from the $k$ dimensional Krylov subspace
$$
\mathcal{K}_k(AA^T,b)=span\{P_k\}.
$$
It is well known from, e.g., \cite{meurant},
that the residual of $y_k^{cgme}$ is
\begin{equation}\label{yksol}
b-AA^T y_k^{cgme}=r_k(AA^T) b,
\end{equation}
where $r_k(\lambda)$ is the $k$-th residual, or Ritz, polynomial with
the normalization $r_k(0)=1$, whose $k$ roots are the Ritz values
$(\bar{\theta}_j^{(k)})^2$ of $AA^T$ with respect to the subspace $span\{P_k\}$;
see \eqref{aat}. Therefore, we have
\begin{equation}\label{rkex}
r_k(\sigma_i^2)=\prod_{j=1}^n \frac{(\bar{\theta}_j^{(k)})^2-\sigma_i^2}
{(\bar{\theta}_j^{(k)})^2},\ i=1,2,\ldots,n.
\end{equation}

From the full SVD \eqref{eqsvd} of $A$, write $U=(U_n,U_{\perp})$. Then
we have $A=U_n\Sigma V^T$, the compact SVD of $A$.
It is straightforward to see that
$$
AA^T(AA^T)^{\dagger}=(AA^T)^{\dagger}AA^T=U_nU_n^T.
$$
Therefore, by $y_{naive}=(AA^T)^{\dagger}b$, premultiplying the two
hand sides of \eqref{yksol} by $(AA^T)^{\dagger}$ yields
\begin{align*}
y_{naive}-U_nU_n^Ty_k^{cgme}&=(AA^T)^{\dagger}r_k(AA^T) b\\
&=r_k(AA^T)(AA^T)^{\dagger} b=r_k(AA^T) y_{naive},
\end{align*}
from which it follows that
\begin{equation}\label{ykex}
U_nU_n^Ty_k^{cgme}=(I-r_k(AA^T))y_{naive}.
\end{equation}
By the SVD \eqref{eqsvd} of $A$, we have
$$
y_{naive}=(AA^T)^{\dagger}b=\sum_{i=1}^n\frac{u_i^Tb}{\sigma_i^2}u_i.
$$
Hence for $ k=1,2,\ldots,n$ from \eqref{rkex} and \eqref{ykex} we obtain
\begin{align}
U_nU_n^T y_k^{cgme}&=\sum\limits_{i=1}^n (1-r_k(\sigma_i^2))\frac{u_i^{T}b}
{\sigma_i^2}u_i \notag\\
&=\sum\limits_{i=1}^nf_i^{(k,cgme)}\frac{u_i^{T}b}{\sigma_i^2}u_i \label{cgmeexp}
\end{align}
with $f_i^{(k,cgme)}$ defined by \eqref{filter}. In terms of
$x_k^{cgme}=A^T y_k^{cgme}$ and $A=U_n\Sigma V^T$, premultiplying
the two hand sides of the above relation by $A^T$ and exploiting
$U_n^TU_n=I$, we have
$$
x_k^{cgme}=A^T y_k^{cgme}=V\Sigma U_n^T y_k^{cgme}
=V\Sigma U_n^T U_n U_n^T y_k^{cgme}=A^T U_nU_n^T y_k^{cgme}.
$$
Then making use of this relation, $A^T u_i=\sigma_i v_i$ and \eqref{cgmeexp},
we obtain \eqref{cgmeexpr}.
\qquad\endproof

Based on Theorems~\ref{interlace}--\ref{theocgme}, we can prove the following
important result.

\begin{theorem}\label{semi}
Let $k_{cgme}^*$ and $k_{lsqr}^*$ be iterations at which
the semi-convergence of CGME and LSQR occurs, respectively, $k_0$
the transition point of the TSVD method. Then
\begin{equation}\label{semicgme}
k_{cgme}^*\leq k_{lsqr}^*\leq k_0,
\end{equation}
that is, the semi-convergence of CGME always occurs no later than that of
LSQR and the TSVD method.
\end{theorem}

{\em Proof}.
The result $k_{lsqr}^*\leq k_0$ has been proved in \cite[Theorem 3.1]{jia18a}.
Next we first prove that $k_{cgme}^*\leq k_0$.

Recall that the best TSVD solution
$$
x_{k_0}^{tsvd}=A_{k_0}^{\dagger}b=\sum_{i=1}^{k_0}\frac{u_i^Tb}{\sigma_i}v_i
$$
and the fact that a 2-norm filtering best possible
solution must capture the $k_0$
dominant SVD components of $A$ and suppress the $n-k_0$ small
SVD components of $A$.

For CGME, from \eqref{thetak} and \eqref{thetak2} we have
$
\bar{\theta}_k^{(k)}<\sigma_k,
$
Therefore, at iteration $k_0+1$ we must have $\bar{\theta}_{k_0+1}^{(k_0+1)}<
\sigma_{k_0+1}$. If the $\bar{\theta}_i^{(k)}$ approximate the large
$\sigma_i$ in natural order for $k=1,2,\ldots,k_0$, then by \eqref{filter}
we have $f_i^{(k,cgme)}\rightarrow 1$ for $i=1,2,\ldots,k$ and
$f_i^{(k,cgme)}\rightarrow 0$ for $i=k+1,\ldots,n$. On the other hand,
by \eqref{filter} we have $f_{k_0+1}^{(k_0+1,cgme)}=\mathcal{O}(1)$.
Compared with the best TSVD solution, by
\eqref{cgmeexpr} the above shows that the CGME iterate $x_k^{cgme}$
captures the $k$ dominant SVD components of $A$ and filters out
the $n-k$ small ones. As a result, $x_k^{cgme}$ improves until iteration $k_0$,
and the semi-convergence of CGME occurs at
iteration $k_{cgme}^*=k_0$.

If the $\bar{\theta}_j^{(k)}$ do not converge to the large singular values of
$A$ in natural order and $\bar{\theta}_k^{(k)}<\sigma_{k_0+1}$
for some iteration $k\leq k_0$ for the first time, then
$x_k^{cgme}$ is already deteriorated by the noise $e$ before iteration $k$:
Suppose that $\sigma_{j^*}<\bar{\theta}_k^{(k)}<\sigma_{k_0+1}$ with
$j^*$ the smallest integer
$j^*>k_0+1$. Then we can easily justify from \eqref{filter} that
$f_i^{(k,cgme)}\in (0,1)$ and tends to zero
monotonically for $i=j^*,j^*+1,\ldots,n$, but
$$
\prod\limits_{j=1}^k\frac{(\bar{\theta}_j^{(k)})^2-\sigma_i^2}
{(\bar{\theta}_j^{(k)})^2}=\frac{(\bar{\theta}_k^{(k)})^2-\sigma_i^2}
{(\bar{\theta}_k^{(k)})^2}\prod\limits_{j=1}^{k-1}
\frac{(\bar{\theta}_j^{(k)})^2-\sigma_i^2}{(\bar{\theta}_j^{(k)})^2}\leq 0,
\ i=k_0+1,\ldots,j^*-1
$$
since the first factor is non-positive and the second factor is positive by
noting that $\bar{\theta}_j^{(k)}>\sigma_i$, $j=1,2,\ldots,k-1$ for
$i=k_0+1,\ldots,j^*-1$. As a result,
$f_i^{(k,cgme)}\geq 1$ for $i=k_0+1,\ldots, j^*-1$, showing
that $x_k^{cgme}$ has been deteriorated by the noise $e$ and
the semi-convergence of CGME has occurred at some iteration $k^*_{cgme}<k_0$.

Finally, we prove $k_{cgme}^*\leq k_{lsqr}^*$.
Notice that $\bar{\theta}_k^{(k)}<\theta_k^{(k)}$ means
that the first iteration $k$ such that $\bar{\theta}_k^{(k)}<\sigma_{k_0+1}$
for CGME is no more than the one such that $\theta_k^{(k)}<\sigma_{k_0+1}$
for LSQR. Therefore,
applying a similar proof to that of the semi-convergence of CGME
to \eqref{eqfilter2}--\eqref{filterlsqr}, it is direct
that the semi-convergence of CGME occurs no later than that of LSQR, i.e.,
$k_{cgme}^*\leq k_{lsqr}^*$.
\qquad\endproof

It is seen from the above proof that, due to $\bar{\theta}_k^{(k)}<\theta_k^{(k)}$,
the semi-convergence of CGME can occur much earlier than that of LSQR.

We can, {\em informally}, deduce more features of CGME.
By definition, the optimality of CGME means that
\begin{equation}\label{cgmelsqr}
\|x_{naive}-x_k^{cgme}\|\leq \|x_{naive}-x_k^{lsqr}\|
\end{equation}
holds unconditionally for $i=1,2,\ldots,n$.
Since $x_k^{cgme}$ and $x_k^{lsqr}$ converge to $x_{true}$ until
iterations $k_{cgme}^*$ and $k_{lsqr}^*$ at which
the semi-convergence of CGME and LSQR occurs, respectively, it is known that,
for $k\leq k_{cgme}^*$ and $k\leq k_{lsqr}^*$, $\|x_{true}-x_k^{cgme}\|$
and $\|x_{true}-x_k^{lsqr}\|$
are negligible relative to $\|x_{naive}-x_{true}\|$, which is supposed very large
in the context of discrete ill-posed problems. As a consequence, we have
\begin{eqnarray}
\|x_{naive}-x_k^{cgme}\|&=&\|x_{naive}-x_{true}+x_{true}-x_k^{cgme}\|
\nonumber \\
&\approx& \|x_{naive}-x_{true}\|+\|x_{true}-x_k^{cgme}\|,\label{naive1}\\
\|x_{naive}-x_k^{lsqr}\|&=&\|x_{naive}-x_{true}+x_{true}-x_k^{lsqr}\| \nonumber\\
&\approx& \|x_{naive}-x_{true}\| +\|x_{true}-x_k^{lsqr}\|.  \label{naive2}
\end{eqnarray}
Since the first terms in the right-hand sides of \eqref{naive1}
and \eqref{naive2} are the same constant, a combination of \eqref{cgmelsqr}
with \eqref{naive1} and \eqref{naive2} means that
\begin{equation}\label{accurcomp}
\|x_{true}-x_k^{cgme}\|\leq \|x_{true}-x_k^{lsqr}\|
\end{equation}
generally holds until $k=k_{cgme}^*$. That is,
$x_k^{cgme}$ should be {\em at least} as accurate as $x_k^{lsqr}$ until
the semi-convergence of CGME occurs. Then for $k>k_{cgme}^*$,
according to Theorem~\ref{semi},
$x_k^{lsqr}$ continues approximating $x_{true}$ as $k$ increases until
iteration $k=k_{lsqr}^*$, at which LSQR ultimately computes
a more accurate approximation $x_{k_{lsqr}^*}^{lsqr}$ to $x_{true}$
than $x_{k_{cgme}^*}^{cgme}$.

We will have more exciting findings. Observe that after
Lanczos bidiagonalization is run $k$ steps, we have
already obtained $\bar{B}_{k+1}$, $P_{k+1}$ and $Q_{k+1}$, but LSQR
and CGME exploit only $B_k, Q_k$ and $\bar{B}_k, Q_k$, respectively.
Since $\alpha_{k+1}>0$ for
$k\leq n-1$, applying the strict interlacing property of singular
values to $B_k$ and $\bar{B}_{k+1}$, we have
\begin{equation}\label{bkbarbk}
\bar{\theta}_1^{(k+1)}>\theta_1^{(k)}>\bar{\theta}_2^{(k+1)}>\cdots>
\bar{\theta}_k^{(k+1)}>\theta_k^{(k)}>\bar{\theta}_{k+1}^{(k+1)}, \
k=1,2,\ldots,n-1.
\end{equation}
Note from \eqref{thetak2} that $\bar{\theta}_i^{(k+1)}<\sigma_i,\
i=1,2,\ldots,k+1$. Combining \eqref{bkbarbk} with \eqref{thetak2},
we see that as approximations
to the first large $k$ singular values $\sigma_i$ of $A$,
although the $k$ singular values $\bar{\theta}_i^{(k)}$
of $\bar{B}_k$ are less accurate than the singular values $\theta_i^{(k)}$
of $B_k$, the first $k$ singular values $\bar{\theta}_i^{(k+1)}$
of $\bar{B}_{k+1}$ are more accurate than the
$\theta_i^{(k)}$ correspondingly.

Based on the above property and \eqref{leftbidiag}, we next show
how to extract a best possible rank $k$ approximation to $A$ from the
available rank $k+1$ matrix $P_{k+1}\bar{B}_{k+1}Q_{k+1}^T=P_{k+1}P_{k+1}^TA$
generated by Algorithm 1.

\begin{theorem}\label{approx}
Let $\bar{C}_k$ be the best rank $k$ approximation to $\bar{B}_{k+1}$ with
respect to the 2-norm. Then for $k=1,2,\ldots,n-1$ we have
\begin{align}
\|A-P_{k+1}\bar{C}_kQ_{k+1}^T\|&\leq \sigma_{k+1}+\gamma_{k+1}^{cgme},\label{lowrank}\\
\|A-P_{k+1}\bar{C}_kQ_{k+1}^T\|&\leq \bar{\theta}_{k+1}^{(k+1)}+\gamma_{k+1}^{cgme},
\label{better}
\end{align}
where $\bar{\theta}_{k+1}^{(k+1)}$ is the smallest singular value of $\bar{B}_{k+1}$
and $\gamma_{k+1}^{cgme}$ is defined by \eqref{cgmeacc}.
\end{theorem}

{\em Proof}. Write $A-P_{k+1}\bar{C}_kQ_{k+1}^T=A-P_{k+1}\bar{B}_{k+1}Q_{k+1}^T+
P_{k+1}(\bar{B}_{k+1}-\bar{C}_k)Q_{k+1}^T$. Then exploiting \eqref{leftbidiag},
we obtain
\begin{align}
\|A-P_{k+1}\bar{C}_kQ_{k+1}^T\| & \leq \|A-P_{k+1}\bar{B}_{k+1}Q_{k+1}^T\|+
\|P_{k+1}(\bar{B}_{k+1}-\bar{C}_k)Q_{k+1}^T\|  \label{barbk}\\
&= \|A-P_{k+1}\bar{B}_{k+1}Q_{k+1}^T\|+
\|P_{k+1}P_{k+1}^TA-P_{k+1}\bar{C}_kQ_{k+1}^T\|. \label{decom2}
\end{align}
By the definition of $C_k$ and \eqref{leftbidiag}, it is easily
justified that $P_{k+1}\bar{C}_kQ_{k+1}^T$ is the
best rank $k$ approximation to $P_{k+1}\bar{B}_{k+1}Q_{k+1}^T=
P_{k+1}P_{k+1}^TA$
in the 2-norm as $P_{k+1}$ and $Q_{k+1}$ are column orthonormal.
Keep in mind that $A_k$ is
the best rank $k$ approximation to $A$. Since $P_{k+1}P_{k+1}^TA_k$ is a rank
$k$ approximation to $P_{k+1}P_{k+1}^TA$, we obtain
\begin{align*}
\|P_{k+1}P_{k+1}^TA-P_{k+1}\bar{C}_kQ_{k+1}^T\| &\leq \|P_{k+1}P_{k+1}^TA
-P_{k+1}P_{k+1}^TA_k\| \\
&=\|P_{k+1}P_{k+1}^T(A-A_k)\|\\
&\leq \|A-A_k\|=\sigma_{k+1}.
\end{align*}
Note that the first term in the right-hand side of \eqref{decom2}
is just $\gamma_{k+1}^{cgme}$. Therefore,
it follow from \eqref{decom2} that \eqref{lowrank} holds.

Since $P_{k+1}$ and $Q_{k+1}$ are column orthonormal and $C_k$ is the best
rank $k$ approximation to $\bar{B}_{k+1}$,
by the orthogonal invariance of 2-norm we obtain
$$
\|P_{k+1}(\bar{B}_{k+1}-\bar{C}_k)Q_{k+1}^T\|=\|\bar{B}_{k+1}-\bar{C}_k\|
=\bar{\theta}_{k+1}^{(k+1)},
$$
which, together with \eqref{barbk}, yields \eqref{better}.
\qquad\endproof

The bound \eqref{better} is always smaller than the bound \eqref{lowrank}
because of $\bar{\theta}_{k+1}^{(k+1)}<\sigma_{k+1}$
from \eqref{thetak} and \eqref{thetak2}. Indeed,
the bound \eqref{lowrank} can be conservative since
we have amplified $\|P_{k+1}(\bar{B}_{k+1}-\bar{C}_k)Q_{k+1}^T\|$
twice and obtained its bound $\sigma_{k+1}$, which might be a considerable
overestimate. Moreover, as we have explained previously,
\eqref{thetak} and \eqref{thetak2} show that
$\bar{\theta}_{k+1}^{(k+1)}>\sigma_n$
may approach $\sigma_n$ for $m=n$ and  $\bar{\theta}_{k+1}^{(k+1)}>0$
can be close to zero arbitrarily for $m>n$.
By definition \eqref{gammak} of $\gamma_k^{lsqr}$, since
$\gamma_{k+1}^{cgme}<\gamma_k^{lsqr}$
(cf. the upper bound of \eqref{cgmelowup}),
$\gamma_k^{lsqr}\geq \sigma_{k+1}>\bar{\theta}_{k+1}^{(k+1)}$ and
$\|A-P_{k+1}\bar{C}_kQ_{k+1}^T\|\geq\sigma_{k+1}$,
the right-hand side of \eqref{better} satisfies
$$
\sigma_{k+1}\leq \bar{\theta}_{k+1}^{(k+1)}+\gamma_{k+1}^{cgme}<2\gamma_k^{lsqr}.
$$
Therefore, $\bar{\theta}_{k+1}^{(k+1)}+\gamma_{k+1}^{cgme}$
is as small as and can even be smaller than
$\gamma_k^{lsqr}$, meaning that $P_{k+1}\bar{C}_kQ_{k+1}^T$ is as
accurate as the rank $k$ approximation $P_{k+1}B_kQ_k^T$ in LSQR.

Define $Q_{n+1}=(Q_n,\mathbf{0})\in \mathbb{R}^{n\times (n+1)}$,
and note from \eqref{bn1} that $\bar{B}_{n+1}=(B_n,\mathbf{0})$.
Recall that the singular values of $\bar{B}_{n+1}$ and $B_n$
are $\bar{\theta}_i^{(n+1)}, \
i=1,2,\ldots,n+1$ and $\theta_i^{(n)},\ i=1,2,\ldots,n$, respectively,
and $\bar{\theta}_i^{(n+1)}=\theta_i^{(n)}=\sigma_i,\ i=1,2,\ldots,n$
and $\bar{\theta}_{n+1}^{(n+1)}=0$. From \eqref{fulllb} and the definition
of $\bar{C}_n$, since $\bar{B}_{n+1}$ is of rank $n$,
we have
$$
\bar{C}_n=\bar{B}_{n+1}
$$
and
$$
A=P_{n+1}B_nQ_n^T=P_{n+1}\bar{B}_{n+1}Q_{n+1}^T=P_{n+1}\bar{C}_nQ_{n+1}^T.
$$
Based on Theorem~\ref{approx} and the analysis followed,
just as done in CGME and LSQR, we can replace $A$ in \eqref{eq1} by the
rank $k$ approximation $P_{k+1}\bar{C}_kQ_{k+1}^T$ and propose
a modified CGME (MCGME) method that solves
\begin{equation}\label{mcgmereg}
\min\|x\| \ \ \mbox{ subject to }\ \ \|P_{k+1}\bar{C}_kQ_{k+1}^T x-b\|=\min
\end{equation}
for the regularized solution $x_k^{mcgme}=Q_{k+1}y_k^{mcgme}$ of \eqref{eq1}
with
\begin{equation}\label{ykmcgme}
y_k^{mcgme}=
\bar{C}_k^{\dagger}P_{k+1}^Tb=\beta_1 \bar{C}_k^{\dagger}e_1^{(k+1)}
\end{equation}
starting with $k=1$ onwards.
MCGME is expected to have the same regularization ability as
LSQR because (i) the $k$ nonzero singular values
$\bar{\theta}_i^{(k+1)}$ of $\bar{C}_k$ are more accurate than the $k$
singular values $\theta_i^{(k)}$ of $B_k$ as approximations to the first
$k$ singular values of $A$ and (ii) $P_{k+1}\bar{C}_kQ_{k+1}^T$ is a rank $k$
approximation which is as accurate as $P_{k+1}B_kQ_k^T$ in LSQR. Regarding
implementations, we comment that the singular values, and left
and right singular vectors of $\bar{C}_k^{\dagger}$ is already available
when $\bar{C}_k$ is extracted from the SVD of $\bar{B}_{k+1}$,
whose computational cost is $\mathcal{O}(k^3)$ flops. As a result,
by \eqref{ykmcgme} we can compute $y_k^{mcgme}$ at cost
of $\mathcal{O}(k^2)$ flops.
A difference from CGME and LSQR is that MCGME seeks $x_k^{mcgme}$
in the $k+1$ dimensional Krylov subspace $\mathcal{K}_{k+1}(A^TA,A^Tb)$
other than in $\mathcal{K}_k(A^TA,A^Tb)$.
Numerical experiments will justify that MCGME
has very comparable regularizing effects to LSQR and can obtain
the best regularized solutions with very similar accuracy to those by
LSQR. We will not consider the by-product MCGME method further in this paper.

$\bar{C}_k$ may have some other potential applications. For example,
when we are required to compute several largest singular triplets of
a large scale matrix $A$, we can use the nonzero singular values of $\bar{C}_k$
to replace the ones of $B_k$ as more accurate approximations to the largest
singular values of $A$ in Lanczos bidiagonaliation type algorithms
\cite{jia03}. In such a way,
exploiting the SVD of $\bar{C}_k$, we can also compute more accurate
approximate left and right singular vectors of $A$ simultaneously.
A development of such modified algorithms is beyond the scope of this paper.

\section{The accuracy of truncated rank $k$ approximate SVDs
by randomized algorithms}\label{randomappro}

In this section, we deviate from the context of Krylov solvers.
Using the analysis approach in the last section, we consider
the accuracy of a truncated rank $k$ SVD approximation to $A$ constructed by
standard randomized algorithms and their improved variants \cite{halko11}.
This topic has been intensively
investigated in recent years; see the survey paper \cite{halko11} and
the references therein. Algorithm 2 is
one of the two basic randomized algorithms from \cite{halko11}
for computing an approximate SVD and extracting a truncated rank $k$ approximate
SVD from it.
A minor difference from the other sections in this paper is that we drop the
restrictions that the singular values of $A$ are simple and $m\geq n$,
that is, the singular values of $A$ are
$\sigma_1\geq \sigma_2\geq\cdots \geq \sigma_{\min\{m,n\}}$.

{\bf Algorithm 2: Randomized approximate SVD of $A$}
\begin{itemize}
\item Input: Given $A\in \mathbb{R}^{m\times n}$, a target rank $k$,
and an oversampling parameter $p$ satisfying $\ell=k+p\leq \min\{m,n\}$.

\item Output: a truncated rank $k$ approximate SVD $A_{(k)}$ of $A$.
\end{itemize}

{\sf Stage A}
\begin{enumerate}
\item  Draw an $n\times \ell$ Gaussian random matrix $\Omega$.

\item Form the $m \times \ell$ matrix $Y = A\Omega$.

\item Compute the compact QR factorization $Y = PR$.
\end{enumerate}

{\sf Stage B}

\begin{enumerate}

\item Form $B=P^TA$.

\item Compute the compact SVD of the $\ell\times n$ matrix $B$:
$B=\widetilde{U}\widetilde{\Sigma}\widetilde{V}^T$.

\item Set $\widehat{U}=P\widetilde{U}$. Compute a rank $\ell$ SVD approximation
$PP^TA=\widehat{U}\widetilde{\Sigma}\widetilde{V}^T$ to $A$.

\item Let $B_{(k)}=\widetilde{U}_k\widetilde{\Sigma}_{(k)}\widetilde{V}_k^T$ be the
best rank $k$ approximation to $B$ with the diagonal $\widetilde{\Sigma}_{(k)}$ being
the first $k$ diagonals of $\widetilde{\Sigma}$, and $\widetilde{U}_k$ and $\widetilde{V}_k$
the first $k$ columns of $\widetilde{U}$ and $\widetilde{V}$, respectively. Form a truncated
rank $k$ SVD approximation $A_{(k)}=PB_{(k)}=
\widehat{U}_k\widetilde{\Sigma}_{(k)}\widetilde{V}_k^T$ to $A$
with $\widehat{U}_k=P\widetilde{U}_k$.
\end{enumerate}

For the approximation accuracy of $A_{(k)}$ to $A$, Halko {\em et al.} \cite{halko11}
establish a fundamental result (cf. Theorem 9.3 there):
\begin{equation}\label{random}
\|A-A_{(k)}\|\leq \sigma_{k+1}+\|(I-PP^T)A\|.
\end{equation}

Assume that the oversampling parameter $p\geq 4$. Making use of
the probability theory, in terms of $\sigma_{k+1}$,
Halko {\em et al.} \cite{halko11} have established a number of bounds
for $\|(I-PP^T)A\|$; see, e.g., Theorems 10.5--10.8 and
Corollary 10.9--10.10 there. However, concerning \eqref{random},
they point out in Remark 9.1
that {\em "In the randomized
setting, the truncation step appears to be less
damaging than the error bound of Theorem 9.3 suggests, but we currently lack a
complete theoretical understanding of its behavior."}
That is to say, the first
term $\sigma_{k+1}$ in \eqref{random} is generally conservative and an overestimate.

Motivated by the proof
of \eqref{better} in Theorem~\ref{approx}, we can improve \eqref{random}
substantially and reveal why \eqref{random} is an overestimate.
Let
\begin{equation}\label{ab}
\widetilde{\sigma}_1\geq\widetilde{\sigma}_2\geq \cdots\geq \widetilde{\sigma}_{k+p}
\end{equation}
be the singular values of $B=P^TA$ defined in Algorithm 2.
It is clear from Algorithm 2 that
$$
P^TAA^TP=BB^T
$$
is an $(k+p)\times (k+p)$ symmetric matrix, which is the projection matrix
of $AA^T$ onto the subspace $span\{P\}$ in the orthonormal basis
of $\{p_i\}_{i=1}^{k+p}$ with $P=(p_1,p_2,\ldots,p_{k+p})$, whose
eigenvalues are $\widetilde{\sigma}_i^2,\ i=1,2,\ldots,k+p$.  Keep in mind
that the eigenvalues of $AA^T$ are $\sigma_i^2,\ i=1,2,\ldots,\min\{m,n\}$
and $m-\min\{m,n\}$ zeros, denoted by
$\sigma_{\min\{m,n\}+1}^2=\cdots=\sigma_m^2=0$ for later use.

\begin{theorem}\label{improve}
For $A\in \mathbb{R}^{m\times n}$, let $P$ and $A_{(k)}$ be defined
as in Algorithm 2, and $\widetilde{\sigma}_{k+1}$ defined as
in \eqref{ab}. Then
\begin{equation}\label{betterbound}
\|A-A_{(k)}\|\leq \widetilde{\sigma}_{k+1}+\|(I-PP^T)A\|
\end{equation}
with
\begin{equation}\label{interrand}
\sigma_{m-p+1}\leq \widetilde{\sigma}_{k+1}\leq \sigma_{k+1}.
\end{equation}
\end{theorem}
{\em Proof.}
Since $P$ is orthonormal, the eigenvalues of $BB^T$ interlace those of $AA^T$ and
satisfy (cf. \cite[p.198, Corollary 4.4]{stewartsun})
$$
\sigma_{m-(k+p)+i}\leq \widetilde{\sigma}_i\leq \sigma_i,\ i=1,2,\ldots,k+p,
$$
from which \eqref{interrand} follows.

From Algorithm 2, we can write
\begin{eqnarray*}
A-A_{(k)}&=&A-PP^TA+PP^TA-A_{(k)} \\
&=&A-PP^TA+
P\widetilde{U}\widetilde{\Sigma}\widetilde{V}^T-P\widetilde{U}_k\widetilde{\Sigma}_{(k)}
\widetilde{V}_{(k)}^T\\
&=& (I-PP^T)A+P(\widetilde{U}\widetilde{\Sigma}\widetilde{V}^T-
\widetilde{U}_k\widetilde{\Sigma}_{(k)}
\widetilde{V}_{(k)}^T).
\end{eqnarray*}
Since $B_{(k)}$ is the best rank $k$ approximation to $B$,
by the column orthonormality of $P$ we obtain
\begin{eqnarray*}
\|A-A_{(k)}\|&\leq &\|(I-PP^T)A\|+\|P(\widetilde{U}\widetilde{\Sigma}\widetilde{V}^T-
\widetilde{U}_k\widetilde{\Sigma}_{(k)}
\widetilde{V}_{(k)}^T)\|\\
&=&\|(I-PP^T)A\|+\|\widetilde{U}\widetilde{\Sigma}\widetilde{V}^T-
\widetilde{U}_k\widetilde{\Sigma}_{(k)}
\widetilde{V}_{(k)}^T\|\\
&=& \|(I-PP^T)A\|+\|B-B_{(k)}\| \\
&=&\|(I-PP^T)A\|+\widetilde{\sigma}_{k+1},
\end{eqnarray*}
which proves \eqref{betterbound}.
\qquad\endproof

\begin{remark}
This theorem indicates that $\widetilde{\sigma}_{k+1}$ never exceeds
$\sigma_{k+1}$ and, for $m, n$ large and $k+p$ small,
it may be much smaller than $\sigma_{k+1}$. Specifically,
$\widetilde{\sigma}_{k+1}$ can be as small as $\sigma_{m-p+1}$.
For $m>n$, whenever $m-p+1>n$, we have $\sigma_{m-p+1}=0$.
Consequently, the
bound \eqref{betterbound} is unconditionally superior to the bound \eqref{random}
and is sharper than the latter when
$\widetilde{\sigma}_{k+1}<\sigma_{k+1}$.
On the other hand, however, note that $\sigma_{k+1}\leq \|A-A_{(k)}\|$.
Therefore, if $\|(I-PP^T)A\|<\sigma_{k+1}$, we must have
$\widetilde{\sigma}_{k+1}\approx\sigma_{k+1}$, that is,
$\widetilde{\sigma}_{k+1}$ dominates the bound \eqref{betterbound}.
Summarizing the above, in response of Remark 9.1 in \cite{halko11},
we conclude that the truncation step does damage the approximation
accuracy of the truncated rank $k$ approximation
when $\|(I-PP^T)A\|<\sigma_{k+1}$ and it is less damaging when
$\|(I-PP^T)A\|\geq \sigma_{k+1}$.
\end{remark}

As we have seen, the column space of $P$ constructed by Algorithm 2 aims to capture
the $(k+p)$-dimensional dominant left singular subspace of $A$.
A variant of it is to capture the $(k+p)$-dimensional right dominant singular
subspace of $A$.
Mathematically, it amounts to applying Algorithm 2 to $A^T$ and
computes a truncated rank $k$ SVD approximation $A_{(k)}$
to $A$ in a similar way. We call such variant Algorithm 3, for which
\eqref{random} now becomes
\begin{equation}\label{randomv}
\|A-A_{(k)}\|\leq \sigma_{k+1}+\|A(I-PP^T)\|
\end{equation}
with the orthonormal $P\in \mathbb{R}^{n\times (k+p)}$.

Note that the eigenvalues of $A^TA$ are $\sigma_i^2,\ i=1,2,\ldots,\min\{m,n\}$
and $n-\min\{m,n\}$ zeros, denoted by $\sigma_{\min\{m,n\}+1}^2=\cdots=\sigma_n^2=0$.
Since the eigenvalues of $(AP)^TAP$ interlace those of $A^TA$, using the
same proof approach as that of Theorem~\ref{improve}, we
can establish the following result.

\begin{theorem}\label{improvev}
For $A\in \mathbb{R}^{m\times n}$, let $P$ and $A_{(k)}$ be defined as
in Algorithm 3, and
$
\widetilde{\sigma}_1\geq\widetilde{\sigma}_2\geq \cdots\geq \widetilde{\sigma}_{k+p}
$
be the singular values of $AP$. Then
\begin{equation}\label{betterboundv}
\|A-A_{(k)}\|\leq \widetilde{\sigma}_{k+1}+\|A(I-PP^T)\|
\end{equation}
with
\begin{equation}\label{interrandv}
\sigma_{n-p+1}\leq \widetilde{\sigma}_{k+1}\leq \sigma_{k+1}.
\end{equation}
\end{theorem}

We comment that, in the case $m<n$, whenever $n-p+1>m$, we have $\sigma_{n-p+1}=0$,
and consequently the bound \eqref{betterboundv} is unconditionally
superior to and can be substantially sharper than the bound \eqref{randomv}
for $m,n$ large and $k+p$ small.

\begin{remark}
If the singular values $\sigma_i$ of $A$ are all simple,
by the strict interlacing properties of eigenvalues,
the singular values of $B$ in Algorithms~2--3 are all simple too,
and the lower and upper bounds in \eqref{interrand}
and \eqref{interrandv} are strict, i.e.,
$\widetilde{\sigma}_{k+1}<\sigma_{k+1}$.
\end{remark}

\begin{remark}
\eqref{random} and \eqref{randomv} and
Theorems~\ref{improve}--\ref{improvev} hold for all the truncated rank
$k$ SVD approximations generated by the enhanced variants of
Algorithm 2--3 in \cite{halko11}, where the unique difference between
the variants is the way that $P$ is generated. More generally,
Theorems~\ref{improve}--\ref{improvev} are true
for arbitrarily given orthonormal $P\in \mathbb{R}^{m\times (k+p)}$
and $P\in \mathbb{R}^{n\times (k+p)}$ with $k+p\leq \min\{m,n\}$, respectively.
\end{remark}


\section{The regularization of LSMR}\label{lsmr}

From Algorithm 1 we obtain
\begin{equation}\label{matrixlsmr}
Q_{k+1}^TA^TAQ_k=(B_k^TB_k,\alpha_{k+1}\beta_{k+1}e_k^{(k)})^T.
\end{equation}
Therefore, from \eqref{yklsmr}, noting
that $Q_{k+1}^TA^Tb=\alpha_1\beta_1 e_1^{(k+1)}$, we have
\begin{equation}\label{lsmrsolution}
x_k^{lsmr}=Q_k(Q_{k+1}^TA^TAQ_k)^{\dagger}Q_{k+1}^TA^Tb,
\end{equation}
which means that LSMR solves the problem
\begin{equation}\label{lsmrrank}
\min\|x\| \ \ \mbox{ subject to }\ \ \|Q_{k+1}Q_{k+1}^TA^TAQ_kQ_k^Tx-A^Tb\|=\min
\end{equation}
for the regularized solution $x_k^{lsmr}$ starting with $k=1$ onwards.
In the meantime,
it is direct to justify that the TSVD solution $x_k^{tsvd}$ solves
the problem
\begin{equation}\label{tsvdreg}
\min\|x\| \ \ \mbox{ subject to }\ \ \|A_k^TA_kx-A^Tb\|=\min
\end{equation}
starting with $k=1$ onwards. Therefore, \eqref{lsmrrank} and
\eqref{tsvdreg} deal with the normal equation $A^TAx=A^T b$ of \eqref{eq1}
by replacing $A^T A$ with its rank $k$ approximations
$Q_{k+1}Q_{k+1}^TA^TAQ_kQ_k^T$ and $A_k^TA_k$, respectively.

In view of \eqref{lsmrrank} and \eqref{tsvdreg},
we need to accurately
estimate the approximation accuracy $\|A^TA-Q_{k+1}Q_{k+1}^TA^TAQ_kQ_k^T\|$ and
investigate how the singular values of $Q_{k+1}^TA^TAQ_k$ approximate
the $k$ large singular values $\sigma_i^2,\ i=1,2,\ldots,k$ of $A^TA$.
We are concerned with some intrinsic relationships between the regularizing
effects of LSMR and those of LSQR and compare the regularization
ability of the two methods.

By \eqref{xk}, \eqref{eqmform1}, \eqref{eqmform2}, \eqref{Bk}
and $P_{k+1}P_{k+1}^Tb=b$, the LSQR iterate
\begin{align*}
x_k^{lsqr}&=Q_kB_k^{\dagger}P_{k+1}^Tb=Q_k(B_k^TB_k)^{-1}B_k^TP_{k+1}^Tb\\
&=Q_k(Q_k^TA^TAQ_k)^{-1}Q_k^TA^TP_{k+1}P_{k+1}^Tb \\
&=Q_k(Q_k^TA^TAQ_k)^{-1}Q_k^TA^Tb,
\end{align*}
which is the solution to the problem
\begin{equation}\label{lsqrrank}
\min\|x\| \ \ \mbox{ subject to }\ \ \|Q_kQ_k^TA^TAQ_kQ_k^T x-A^Tb\|=\min
\end{equation}
that replaces $A^TA$ by its rank $k$ approximation $Q_kQ_k^TA^TAQ_kQ_k^T
=Q_kB_k^TB_kQ_k^T$ in the normal equation $A^T Ax=A^T b$.
In this sense, the accuracy of such rank $k$ approximation
is measured in terms of $\|A^TA-Q_kQ_k^TA^TAQ_kQ_k^T\|$ for LSQR.

Firstly, we present the following result, which compares the
accuracy of two rank $k$ approximations involved in LSMR and LSQR
in the sense of solving the normal equation $A^TA x=A^T b$.

\begin{theorem}\label{lsqrmr}
For the rank $k$ approximations to $A^TA$ in \eqref{lsmrrank}
and \eqref{lsqrrank}, $k=1,2,\ldots,n-1$, we have
\begin{align}\label{lsqrmrest}
\|A^TA-Q_{k+1}Q_{k+1}^TA^TAQ_kQ_k^T\|& <
\|A^TA-Q_kQ_k^TA^TAQ_kQ_k^T\|.
\end{align}
\end{theorem}

{\em Proof.}
For the orthogonal matrix $Q_n$ generated by Algorithm 1, noticing that
$\alpha_{n+1}=0$, from \eqref{eqmform1} and \eqref{eqmform2} we obtain
$Q_n^TA^TAQ_n=B_n^TB_n$ and
\begin{align}
\|A^TA-Q_{k+1}Q_{k+1}^TA^TAQ_kQ_k^T\|&=
\|Q_n^T(A^TA-Q_{k+1}Q_{k+1}^TA^TAQ_kQ_k^T)Q_n\| \nonumber\\
&=\|B_n^TB_n-(I,\mathbf{0})^T(B_k^TB_k,\alpha_{k+1}
\beta_{k+1}e_k)^T(I,\mathbf{0})\| \nonumber \\
&=\|F_k\|, \ k=1,2,\ldots,n-1, \label{fklsmr}
\end{align}
where
\begin{align}\label{fk}
F_k&=\left(\begin{array}{ccccc}
\alpha_{k+1}\beta_{k+1} & & & &\\
\alpha_{k+1}^2+\beta_{k+2}^2 &\alpha_{k+2}\beta_{k+2} &&&\\
\alpha_{k+2}\beta_{k+2} &\alpha_{k+2}^2+\beta_{k+3}^2&\ddots & &\\
& \alpha_{k+3}\beta_{k+3} & \ddots& & \\
& & &\alpha_{n-1}\beta_{n-1}& \\
 & & \ddots& \alpha_{n-1}^2+\beta_n^2 &\alpha_n\beta_n\\
  & & &\alpha_n\beta_n&\alpha_n^2+\beta_{n+1}^2\\
\end{array}
\right)\\
&=\left(\begin{array}{c}
\alpha_{k+1}\beta_{k+1}(e_1^{(n-k)})^T\\
G_k^TG_k
\end{array}
\right)\in \mathbb{R}^{(n-k+1)\times (n-k)} \label{fkgk}
\end{align}
is the matrix by deleting
the $(k+1)\times k$ leading principal matrix of the symmetric tridiagonal
matrix $B_n^TB_n$ and the first $k-1$ zero rows and $k$ zero columns
of the resulting matrix, where $G_k$ is defined by \eqref{gk1}
and $e_1^{(n-k)}$ are the first canonical vector of $\mathbb{R}^{n-k}$.

On the other hand, it is direct to verify that
\begin{align}
\|A^TA-Q_kQ_k^TA^TAQ_kQ_k^T\|&=
\|Q_n^T(A^TA-Q_kQ_k^TA^TAQ_kQ_k^T)Q_n\| \nonumber\\
&=\|B_n^TB_n-(I,\mathbf{0})^TB_k^TB_k(I,\mathbf{0})\| \nonumber\\
&=\|F_k^{\prime}\|, \label{fklsqr}
\end{align}
where $F_k^{\prime}=\left(\alpha_{k+1}\beta_{k+1}
e_2^{(n-k+1)},F_k\right)\in \mathbb{R}^{(n-k+1)\times (n-k+1)}$
with $e_2^{(n-k+1)}$ being the second canonical vector of $\mathbb{R}^{n-k+1}$.

From \eqref{fk} and \eqref{fkgk}, we obtain
\begin{align}
F_k^{\prime} (F_k^{\prime})^T&=
(\alpha_{k+1}\beta_{k+1} e_2^{(n-k+1)},F_k) (\alpha_{k+1}\beta_{k+1} e_2^{(n-k+1)},
F_k)^T \nonumber\\
&= F_kF_k^T+\alpha_{k+1}^2\beta_{k+1}^2 e_2^{(n-k+1)}(e_2^{(n-k+1)})^T.
\label{fkprime}
\end{align}
Since $G_k^TG_k$ is unreduced symmetric tridiagonal, its eigenvalues
are all simple. Observe from \eqref{fkgk} that
\begin{equation}\label{fkfk}
F_k^TF_k=(G_k^TG_k)^2+\alpha_{k+1}^2\beta_{k+1}^2 e_1^{(n-k)} (e_1^{(n-k)})^T,
\ k=1,2,\ldots,n-1.
\end{equation}
Therefore, we know from \cite[p.218]{demmel} that the eigenvalues of $F_k^TF_k$
strictly interlace those of $(G_k^T G_k)^2$ and are all simple. Furthermore,
we see from \eqref{gk1} that $G_k$ is of full column rank,
which means that the eigenvalues of $F_k^TF_k$ are all
{\em positive}.

Note that the eigenvalues of $F_kF_k^T$ are those of $F_k^TF_k$ and zero.
As a result, the eigenvalues of $F_kF_k^T$ are all simple.
According to \cite[p.218]{demmel},  we know from \eqref{fkprime} that
the eigenvalues of $F_k^{\prime} (F_k^{\prime})^T$ strictly interlace those
of $F_kF_k^T$. Therefore, we obtain
$$
\|F_k^{\prime}\|^2=\|F_k^{\prime}(F_k^{\prime})^T\|>\|F_kF_k^T\|=\|F_k\|^2,
$$
which, from \eqref{fklsmr} and \eqref{fklsqr}, establishes \eqref{lsqrmrest}.
\qquad\endproof

This theorem indicates that, as far as solving $A^TAx=A^Tb$
is concerned, the rank $k$ approximation in LSMR is more accurate than
that in LSQR.

Recall that \eqref{gammak} measures the quality of the rank $k$
approximation involved in LSQR for the regularization
problem \eqref{lsqrreg}. We now
estimate the approximation accuracy of $Q_{k+1}Q_{k+1}^TA^TAQ_kQ_k^T$
to $A^TA$ in terms of $(\gamma_k^{lsqr})^2$.

\begin{theorem}\label{aprod}
For $k=1,2,3,\ldots,n-1$, let $\gamma_k^{lsqr}$ be defined as \eqref{gammak}.
For $k=2,3,\ldots,n-1$ we have
\begin{equation}\label{aproderror}
(\gamma_k^{lsqr})^2< \|A^TA-Q_{k+1}Q_{k+1}^TA^TAQ_kQ_k^T\|\leq
\sqrt{1+m_k(\gamma_{k-1}^{lsqr}/\gamma_k^{lsqr})^2}(\gamma_k^{lsqr})^2
\end{equation}
with $0< m_k<1$ and $\gamma_0^{lsqr}=\|A\|$.
For $k=1,2,\ldots,n-2$, the strict monotonic decreasing property holds:
\begin{equation}\label{monolsmr}
\|A^TA-Q_{k+1}Q_{k+1}^TA^TAQ_kQ_k^T\|<
\|A^TA-Q_{k+2}Q_{k+2}^TA^TAQ_{k+1}Q_{k+1}^T\|.
\end{equation}
\end{theorem}

{\em Proof}.
Combining \eqref{fkgk} with \eqref{gk} and \eqref{alphagamma},
for $k=2,3,\ldots,n-1$ we obtain from \cite[p.98]{wilkinson}
and \cite[p.218]{demmel} that
\begin{equation}\label{fgk}
\|F_k\|^2=\|G_k\|^4+m^{\prime}_k\alpha_{k+1}^2\beta_{k+1}^2\leq
(\gamma_k^{lsqr})^4+m_k(\gamma_{k-1}^{lsqr}\gamma_k^{lsqr})^2
\end{equation}
with $0< m^{\prime}_k\leq 1$ and $0< m_k<m^{\prime}_k$,
from which the lower and upper bounds
in \eqref{aproderror} follow directly.

For $k=1$, the equality in \eqref{fgk} is still true.
From \eqref{alphagamma}, we have $\alpha_2<\gamma_1^{lsqr},\
\beta_2<\|A\|=\gamma_0^{lsqr}$. Therefore,
we obtain
$$
(\gamma_1^{lsqr})^4<\|F_1\|^2=\|G_1\|^4+m_1^{\prime} \alpha_2^2\beta_2^2
\leq (\gamma_1^{lsqr})^4+m_1(\gamma_0^{lsqr}\gamma_1^{lsqr})^2,
$$
from which it follows that \eqref{aproderror} holds for $k=1$.

From \eqref{fk}, we see that
$F_{k+1}$ is the matrix that first deletes
the first column and row of $F_k$ and then deletes the first zero column and row of the
resulting matrix. Therefore, applying the interlacing property of singular values
to $F_{k+1}$ and $F_k$ yields
$$
\|F_k\|\leq \|F_{k+1}\|.
$$

We next prove that the above "$\leq$" is the strict "$<$".
Since $B_n^TB_n=Q_n^TA^TAQ_n$ is an unreduced symmetric tridiagonal matrix,
its singular values $\sigma_i^2,\ i=1,2,\ldots,n$ are simple.
Observe that $F_k$ is the matrix deleting the first $k$ columns of $B_n^TB_n$ and
the first $k$ zero rows of the resulting matrix. Consequently,
the singular values $\zeta_i^{(k)},\, i=1,2,\ldots,n-k$
of $F_k$ strictly interlace the simple singular values
$\sigma_i^2,\ i=1,2,\ldots,n$ of $B_n^TB_n$ and are thus simple
for $k=1,2,\ldots,n-1$. Moreover, the singular values of
$F_{k+1}$ strictly interlace those of $F_k$, which means that
$\zeta_1^{(k)}<\zeta_1^{(k+1)}$, i.e., $\|F_k\|<\|F_{k+1}\|$, which proves
\eqref{monolsmr}.
\qquad\endproof

\begin{remark}
According to the results and analysis in \cite{jia18b},
we have $\gamma_{k-1}^{lsqr}/\gamma_k^{lsqr}\sim \rho$ for severely ill-posed
problems, and $\gamma_{k-1}^{lsqr}/\gamma_k^{lsqr}\sim (k/(k-1))^{\alpha}$ for
moderately and mildly ill-posed problems. Therefore,
the lower and upper bounds of \eqref{aproderror} indicate
that $\|A^T A-Q_{k+1}Q_{k+1}^TA^TAQ_kQ_k^T\|\sim (\gamma_k^{lsqr})^2$.
\end{remark}

Finally, let us investigate the relationship between the singular values of
rank $k$ approximation matrices in LSMR and LSQR. From
\eqref{matrixlsmr} and \eqref{Bk}, we know that they are
the singular values of $(B_k^TB_k,\alpha_{k+1}\beta_{k+1}e_k^{(k)})^T$ and $B_k^TB_k$,
respectively.

\begin{theorem}\label{lsmrlsqr}
Let $(\widetilde{\theta}_1^{(k)})^2> (\widetilde{\theta}_2^{(k)})^2>
\cdots>(\widetilde{\theta}_k^{(k)})^2$ be the singular values of
$(B_k^TB_k,\alpha_{k+1}\beta_{k+1}e_k^{(k)})^T$.
Then for $i=1,2,\ldots,k$ we have
\begin{align}
\theta_i^{(k)}&<\widetilde{\theta}_i^{(k)}<\sigma_i, \label{lsmrqr}\\
(\widetilde{\theta}_i^{(k)})^2&<(\theta_i^{(k)})^2+\gamma_k^{lsqr}
\gamma_{k-1}^{lsqr}. \label{lsmrvalue}
\end{align}
\end{theorem}

{\em Proof.}
Observe that $(B_k^TB_k,\alpha_{k+1}\beta_{k+1}e_k^{(k)})^T$ is the matrix
consisting of the first $k$ columns of $B_n^TB_n$ and deleting
the last $n-k-1$ zero rows of the resulting
matrix. As a result, since $\sigma_i,\ i=1,2,\ldots,n$, are simple,
the singular values $(\widetilde{\theta}_i^{(k)})^2$ of
$(B_k^TB_k,\alpha_{k+1}\beta_{k+1}e_k^{(k)})^T$ strictly
interlace the singular values $\sigma_i^2$ of $B_n^TB_n$:
$$
\sigma_{n-k+i}^2< (\widetilde{\theta}_i^{(k)})^2 <\sigma_i^2,\
i=1,2,\ldots,k
$$
and are simple, which means the upper bound \eqref{lsmrqr}.

Note that
$(B_k^TB_k,\alpha_{k+1}\beta_{k+1}e_k^{(k)})^T(B_k^TB_k,\alpha_{k+1}\beta_{k+1}e_k^{(k)})$
has the $k+1$ eigenvalues $(\widetilde{\theta}_i^{(k)})^4$ and zero, and $(B_k^TB_k)^T(B_k^TB_k)=
(B_k^TB_k)^2$ is its
$k\times k$ leading principal submatrix and has $k$ simple eigenvalues $(\theta_i^{(k)})^4$.
Therefore, $(\theta_i^{(k)})^4$ strictly interlace  $(\widetilde{\theta}_i^{(k)})^4$ and zero,
which proves the lower bound of \eqref{lsmrqr}.

On the other hand, we have
$$
(B_k^TB_k,\alpha_{k+1}\beta_{k+1}e_k^{(k)})(B_k^TB_k,\alpha_{k+1}
\beta_{k+1}e_k^{(k)})^T
= (B_k^TB_k)^2+\alpha_{k+1}^2\beta_{k+1}^2e_k^{(k)}(e_k^{(k)})^T.
$$
Recall \eqref{gk1} that $\alpha_{k+1}<\gamma_k^{lsqr}$
and $\beta_{k+1}<\gamma_{k-1}^{lsqr}$. By standard perturbation theory, we obtain
$$
(\widetilde{\theta}_i^{(k)})^4-(\theta_i^{(k)})^4\leq \alpha_{k+1}^2\beta_{k+1}^2
<(\gamma_k^{lsqr}\gamma_{k-1}^{lsqr})^2, \  i=1,2,\ldots,k,
$$
from which it follows that \eqref{lsmrvalue} holds.
\qquad\endproof

\begin{remark}\label{semilsmr}
\eqref{lsmrqr} indicates that $\widetilde{\theta}_i^{(k)},\ 1=1,2,\ldots,k$
approximate the first $k$ large singular values $\sigma_i$ more accurately
than $\theta_i^{(k)}$. Particularly, since
$\theta_k^{(k)}<\widetilde{\theta}_k^{(k)}$, the first iteration step
$k$ such that $\widetilde{\theta}_{k}^{(k)}<\sigma_{k_0+1}$
must be no smaller than the $k$ such that $\theta_{k}^{(k)}<\sigma_{k_0+1}$.
A combination of this and the previous analysis
on the semi-convergence of CGME and LSQR implies that
the semi-convergence of LSMR must occur no sooner than that of LSQR.
On the other hand, \eqref{lsmrvalue} shows that $\widetilde{\theta}_i^{(k)}$
is bounded from the above by $\theta_i^{(k)}$
as an approximation to $\sigma_i$, which and \eqref{lsmrqr} imply
that $\widetilde{\theta}_i^{(k)}$ and $\theta_i^{(k)}$ interact and
$\theta_i^{(k)}$ cannot be considerably more accurate than
$\widetilde{\theta}_i^{(k)}$ as approximations to the large
singular values of $A$ for $i=1,2,\ldots,k$.
\end{remark}

\begin{remark}
A combination of Theorem~\ref{lsqrmr} and the above two remarks means that the
regularizing effects of LSMR are not inferior to
those of LSQR and the best regularized solutions by LSMR are at least as
accurate as those by LSQR, that is, LSMR has the same regularization ability
as that of LSQR. Particularly,
from the results on LSQR in Section \ref{lsqr}, we conclude that LSMR has
the full regularization for severely or moderately ill-posed problems with
suitable $\rho>1$ or $\alpha>1$.
\end{remark}

A final note is that Huang and Jia \cite{huangjia17} have derived
the eigendecomposition, i.e., equivalent SVD,
filtered expansion of MINRES iterates for $Ax=b$
with $A$ symmetric; see Theorem~3.1 there. The result can be directly
adapted to the LSMR iterates $x_k^{lsmr}$ by keeping in mind that LSMR
is mathematically equivalent to MINRES applied to the specific
symmetric positive definite linear system $A^TAx=A^Tb$.

\section{Numerical experiments}\label{numer}

All the computations are carried out in Matlab R2017b on the
Intel Core i7-4790k with CPU 4.00 GHz processor and 16 GB RAM
with the machine precision
$\epsilon_{\rm mach}= 2.22\times10^{-16}$ under the Miscrosoft
Windows 8 64-bit system.

We have tested LSQR, CGME, LSMR and MCGME on almost all the 1D and 2D
problems from \cite{berisha,hansen07,hansen12} and have
observed similar phenomena.
For the sake of length, we list only some of
them in Table~\ref{tab1}, where each problem takes its default parameter(s).
We mention that the relatively easy 1D problems are all
from \cite{hansen07,hansen12}, where
{\sf shaw}, {\sf gravity} and {\sf baart} are severely
ill-posed and {\sf phillips}, {\sf heat} and and {\sf deriv2}
are moderately. The 2D image deblurring problems {\sf blur}, {\sf fanbeamtomo}
and {\sf seismictomo} are also from \cite{hansen07,hansen12}, and the
other 2D problems are from \cite{berisha}.
We notice that for {\sf blur}, {\sf fanbeamtomo},
although the orders $m$ and $n$ are already tens of thousands,
their condition numbers $\sigma_1/\sigma_n$ are only 31.5 and 2472,
respectively, which, intuitively, do not satisfy the definition of a
discrete ill-posed problem whose singular values decay and are centered
at zero, so that the ratio $\sigma_1/\sigma_n$ is very large.
For each test problem, we compute $b_{true}=Ax_{true}$ and add a Gaussian
white noise $e$ with zero mean to $b_{true}$ by prescribing the relative noise
level
\begin{equation}\label{noiselevel}
\varepsilon=\frac{\|e\|}{\|b_{true}\|}.
\end{equation}

\graphicspath{{figurecgme/}}

\begin{table}[h]
    \centering
    \caption{The description of test problems.}
    \begin{tabular}{lll}
     \hline
     Problem        & Description                                & Size of $m,\ n$ \\
     \hline
     {\sf shaw}   &1D image restoration model & $m=n=5000$\\
     {\sf gravity} &1D gravity surveying problem & $m=n=5000$\\
     {\sf baart} & 1D image deblurring & $m=n=5000$\\
     {\sf phillips} &phillips' famous test problem & $m=n=5000$ \\
     {\sf heat} & Inverse heat problem  & $m=n=5000 $\\
     {\sf deriv2} &Computation of second derivative &$m=n=10000$\\
     {\sf AtmosphericBlur10}     & Spatially Invariant Gaussian Blur   &$m=n=65536$\\
     {\sf AtmosphericBlur30}     & Spatially Invariant Gaussian Blur   &$m=n=65536$\\
      {\sf GaussianBlur420}     & Spatially Invariant Atmospheric   &$m=n=65536$\\
      &  Turbulence Blur&\\
      {\sf GaussianBlur422}     & Spatially Invariant Atmospheric   &$m=n=65536$\\
      &  Turbulence Blur& \\
{\sf VariantGaussianBlur1}   &Spatially Variant Gaussian Blur  & $m=n=99856$ \\
{\sf VariantGaussianBlur2}   &Spatially Variant Gaussian Blur  & $m=n=99856$ \\
{\sf VariantMotionBlur\_large}   &Spatially Variant Motion Blur  & $m=n=65536$ \\
{\sf VariantMotionBlur\_medium}   &Spatially Variant Motion Blur & $m=n=65536$ \\
{\sf blur}   &2D image restoration & $m=n=22500$ \\
{\sf fanbeamtomo}   &2D fan-beam tomography problem  & $61200\times 14400$ \\
{\sf seismictomo} & 2D seismic tomography & $20000\times 10000$\\
     \hline
   \end{tabular}
   \label{tab1}
\end{table}

We use the code {\sf lsqr\_b.m} of \cite{hansen07}, where the reorthogonalization
is exploited during Lanczos bidiagonalization
in order to maintain the numerical orthogonality of $P_{k+1}$ and
$Q_k$. We have written the Matlab codes of CGME, LSMR and MCGME based on the
same Lanczos bidiagonalization process used in {\sf lsqr\_b.m}.

For all the 1D problems and the 2D {\sf seismictomo}, we report the results
on them for $\varepsilon=10^{-3}$; for all the 2D problems except
{\sf blur} and {\sf fanbeamtomo}, we report the results on
them for $\varepsilon=5\times 10^{-3}$. For several other
$\varepsilon\in [10^{-3},5\times 10^{-2}]$, we have the same findings.
For {\sf blur} and {\sf fanbeamtomo}, however, we will observe some
fundamental distinctions between the convergence features
for $\varepsilon$ lying in this practical interval.
Figures~\ref{fig1}--\ref{fig6} depict the convergence processes of
LSQR, CGME, LSMR and MCGME, and we give
some key details, including the iterations $k^*$ at which the
semi-convergence of an algorithm occurs and the relative error
of the best regularized solution obtained by each algorithm, which is
defined by
$$
\frac{\|x_{k^*}^{lsqr}-x_{true}\|}{\|x_{true}\|}
$$
for LSQR. Similar relative errors are defined for CGME, LSMR and MCGME with
the superscript ``$lsqr$'' replaced by ``$cgme$'', ``$lsmr$'' and ``$mcgme$'',
respectively. In addition, as a comparison standard on the solution accuracy,
we depict the semi-convergence process of the TSVD method for {\sf blur} and
{\sf seismictomo}, and report the relative errors of the best TSVD regularized
solutions $x_{k_0}^{tsvd}$ with $k_0$ the transition point at which the semi-convergence
of TSVD occurs. For the other nine larger 2D problems, we cannot compute
the SVDs of the matrices due to out of memory in our computer.
We mention that for the first six 1D test problems we have found
that the best regularized solutions obtained by TSVD method have the same accuracy
as those by LSQR, where the $k_0$ are very small relative to $n$ and all the
$k^*\leq k_0$ correspondingly. We omit the results on the 1D problems
obtained by the TSVD method.

\begin{figure}
\begin{minipage}{0.48\linewidth}
  \centerline{\includegraphics[width=6.0cm,height=4.5cm]{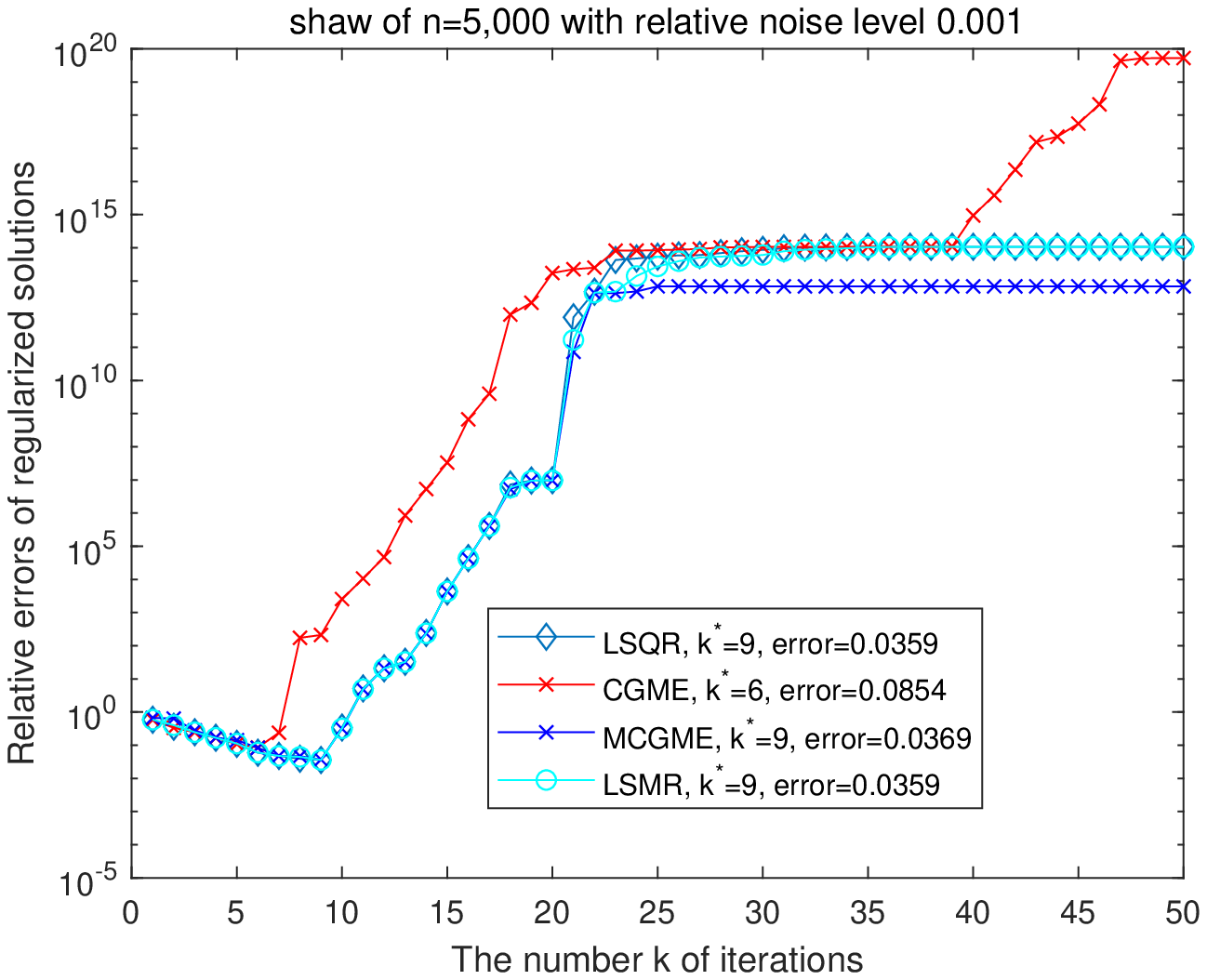}}
  \centerline{(a)}
\end{minipage}
\hfill
\begin{minipage}{0.48\linewidth}
  \centerline{\includegraphics[width=6.0cm,height=4.5cm]{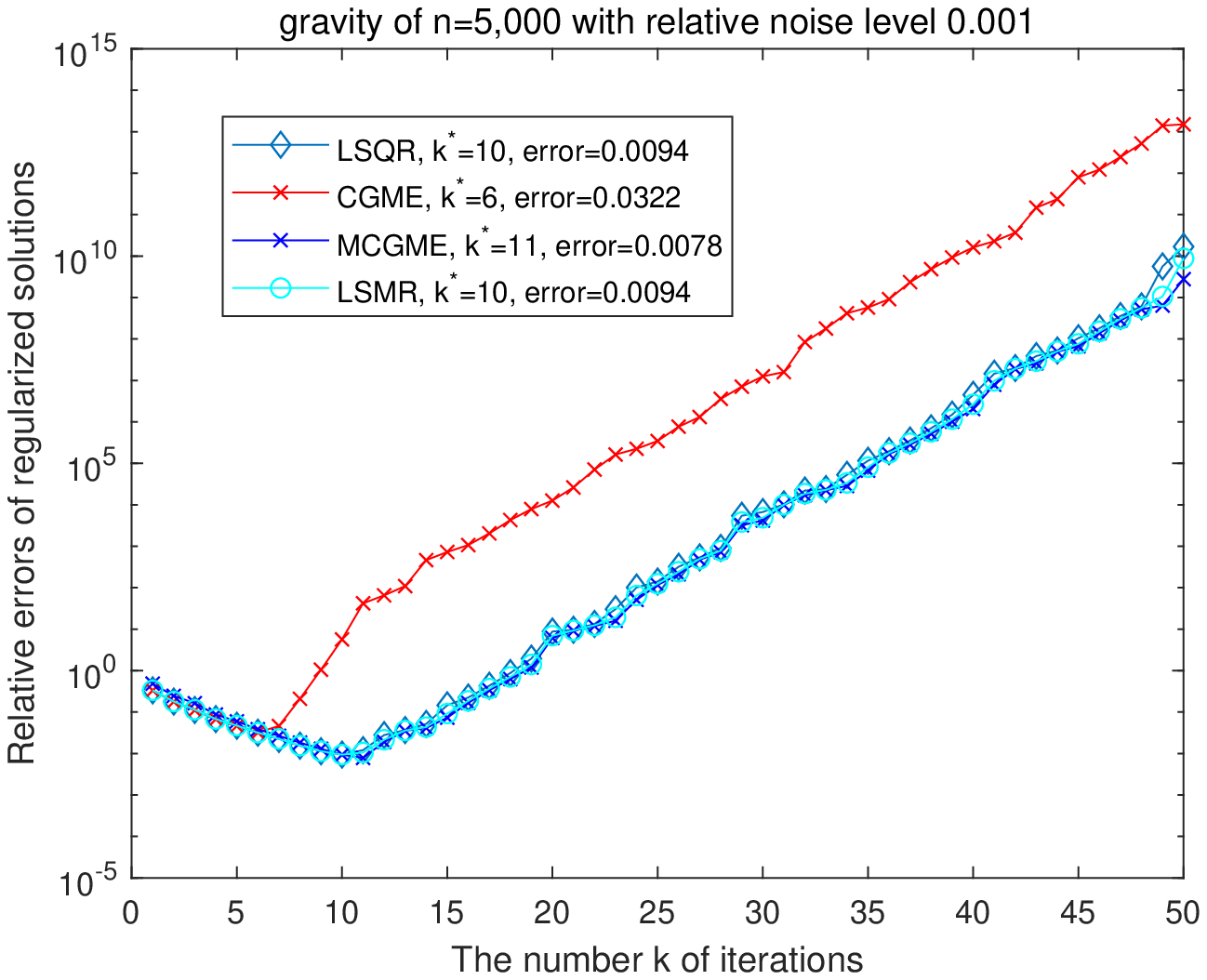}}
  \centerline{(b)}
\end{minipage}
\vfill
\begin{minipage}{0.48\linewidth}
  \centerline{\includegraphics[width=6.0cm,height=4.5cm]{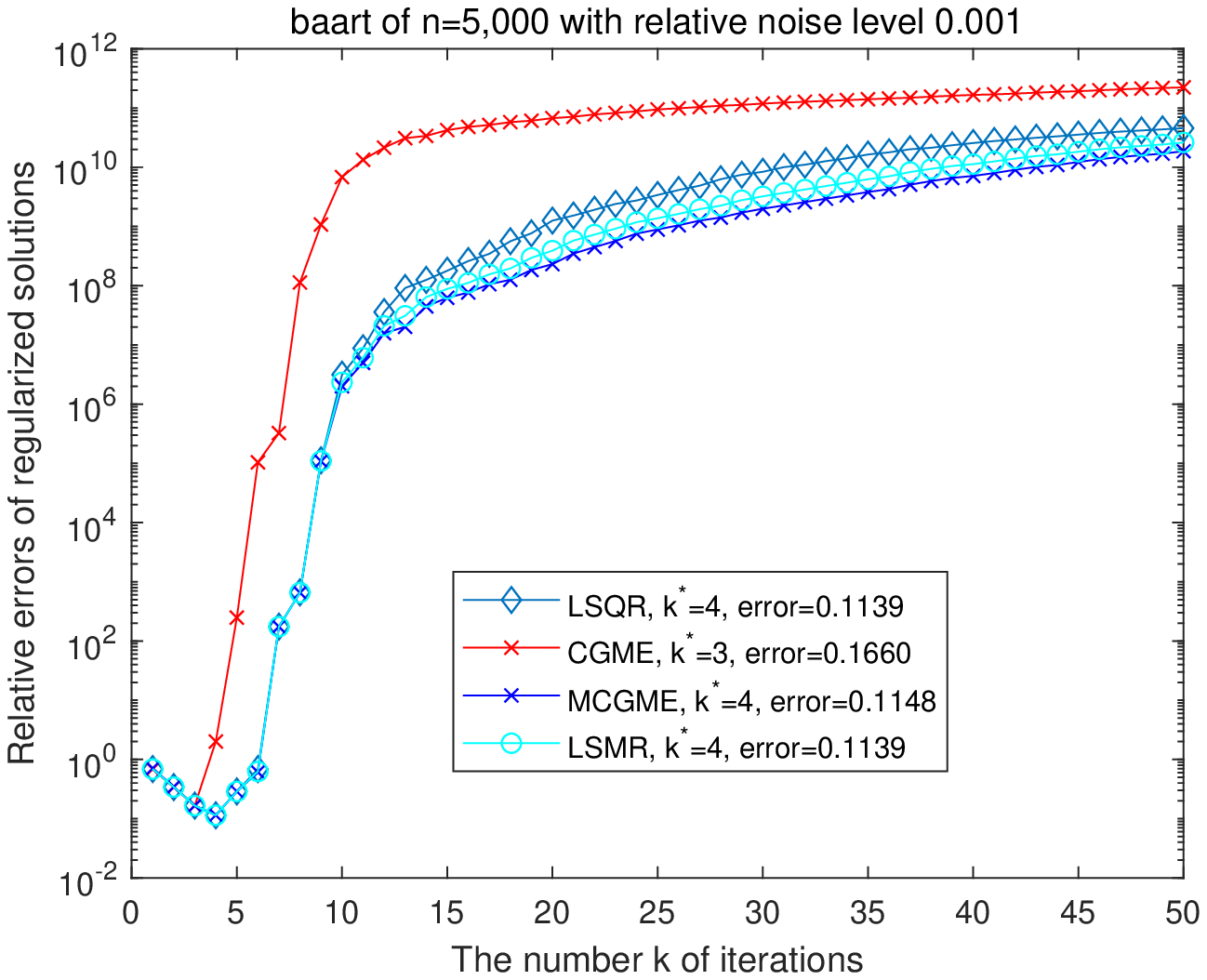}}
  \centerline{(c)}
\end{minipage}
\hfill
\begin{minipage}{0.48\linewidth}
  \centerline{\includegraphics[width=6.0cm,height=4.5cm]{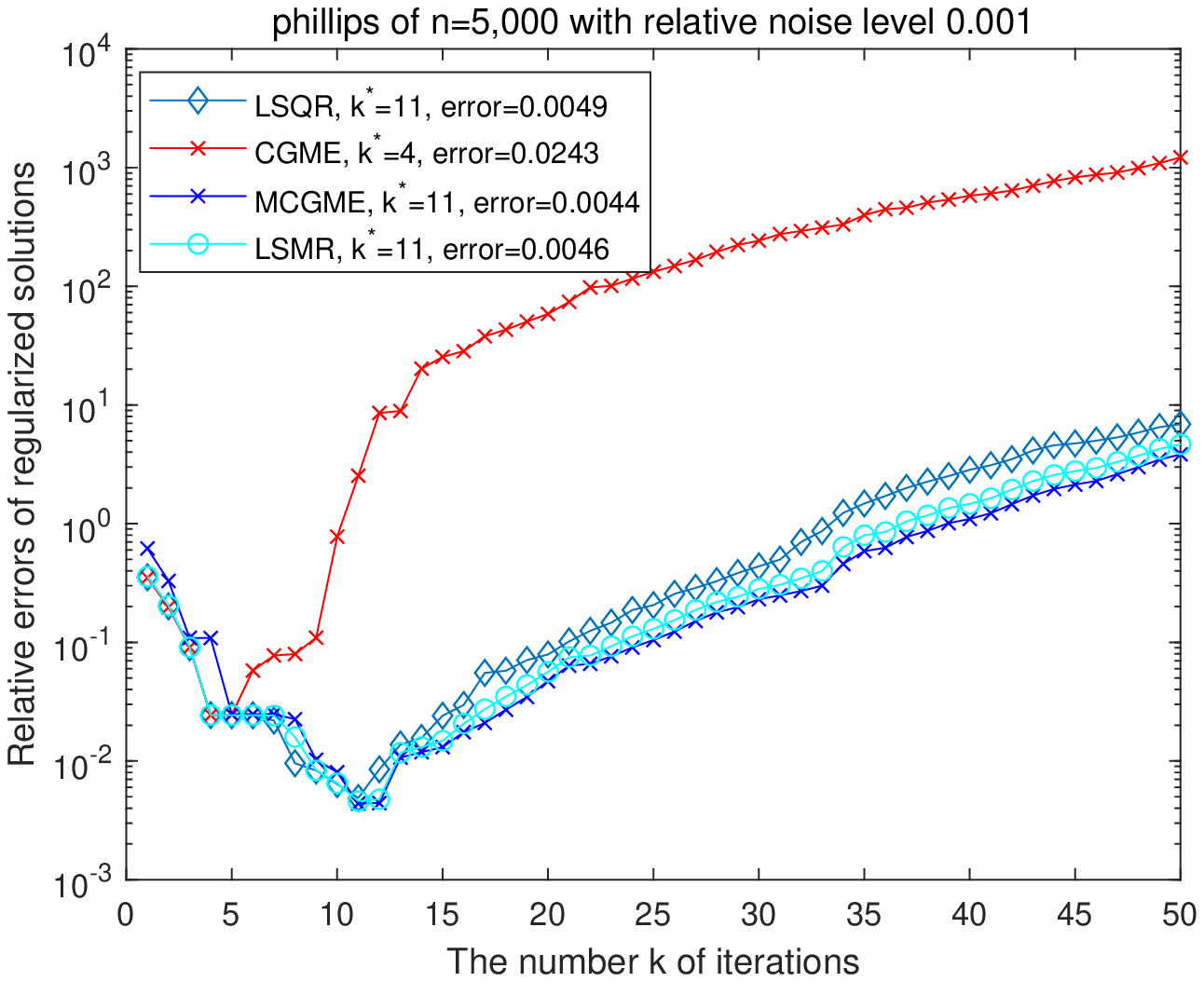}}
  \centerline{(d)}
\end{minipage}
\vfill
\begin{minipage}{0.48\linewidth}
  \centerline{\includegraphics[width=6.0cm,height=4.5cm]{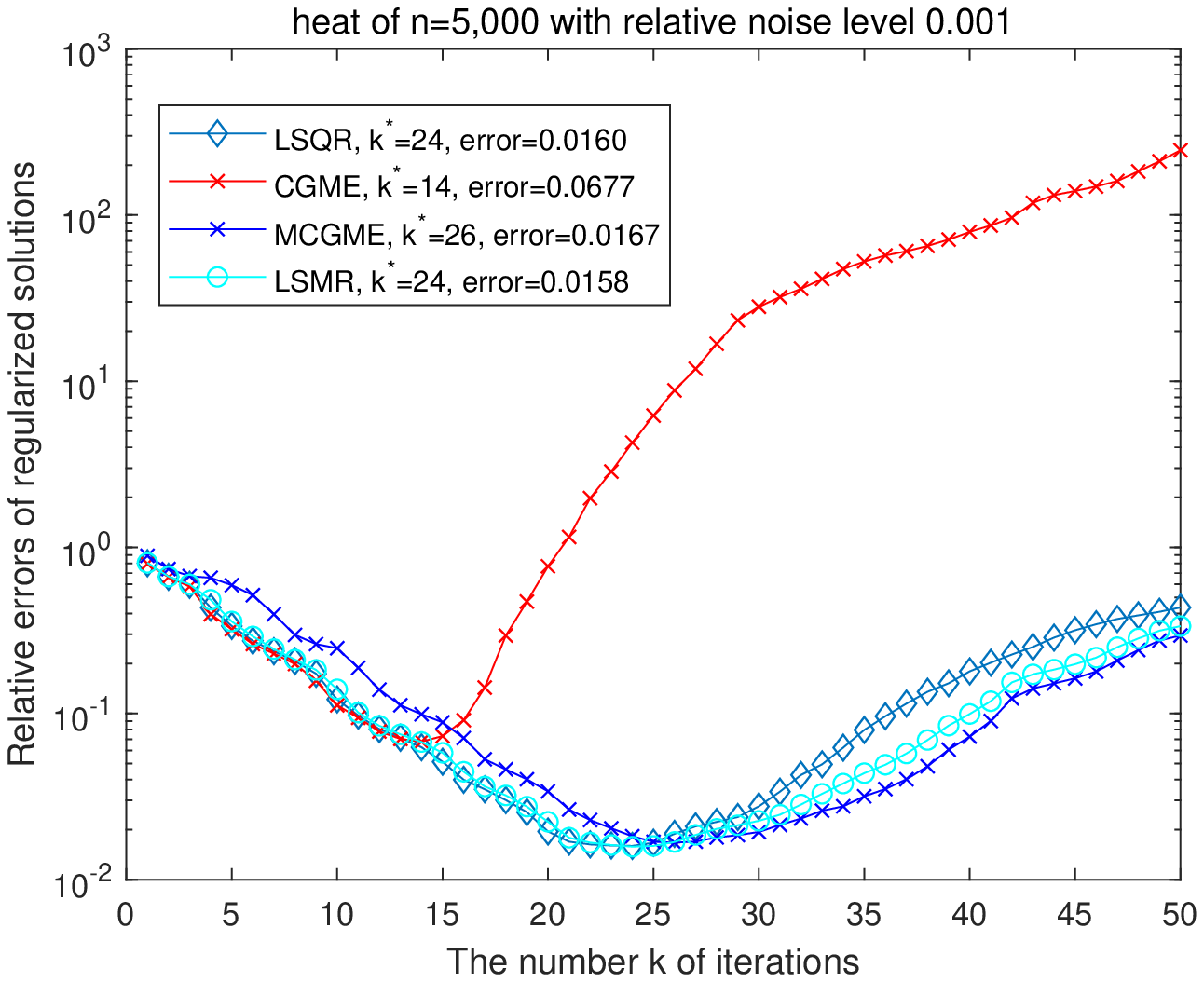}}
  \centerline{(e)}
\end{minipage}
\hfill
\begin{minipage}{0.48\linewidth}
  \centerline{\includegraphics[width=6.0cm,height=4.5cm]{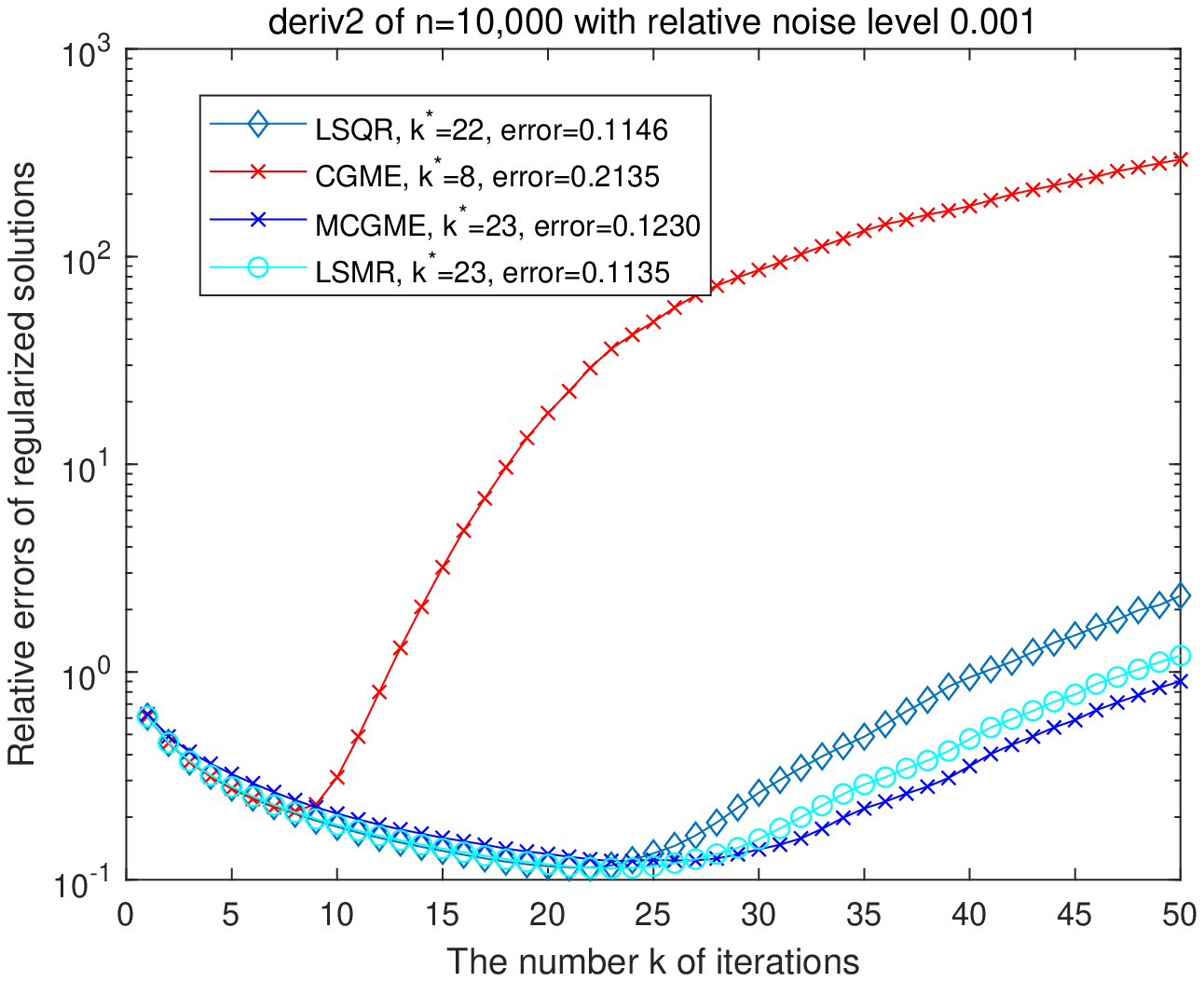}}
  \centerline{(f)}
\end{minipage}
\caption{1D problems with the relative noise level $\varepsilon=10^{-3}$.} \label{fig7}
\end{figure}

\begin{figure}
\begin{minipage}{0.48\linewidth}
  \centerline{\includegraphics[width=6.0cm,height=4.5cm]{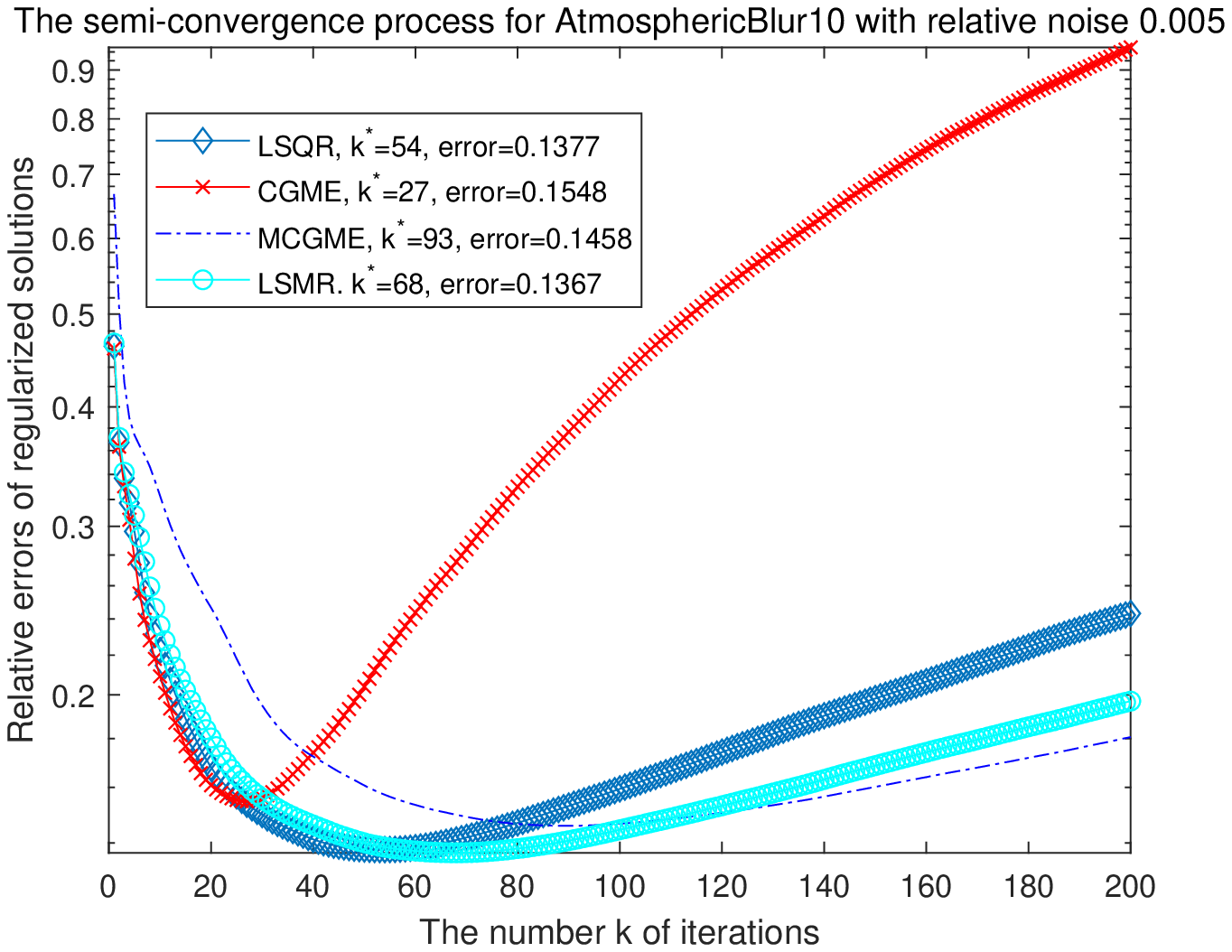}}
  \centerline{(a)}
\end{minipage}
\hfill
\begin{minipage}{0.48\linewidth}
  \centerline{\includegraphics[width=6.0cm,height=4.5cm]{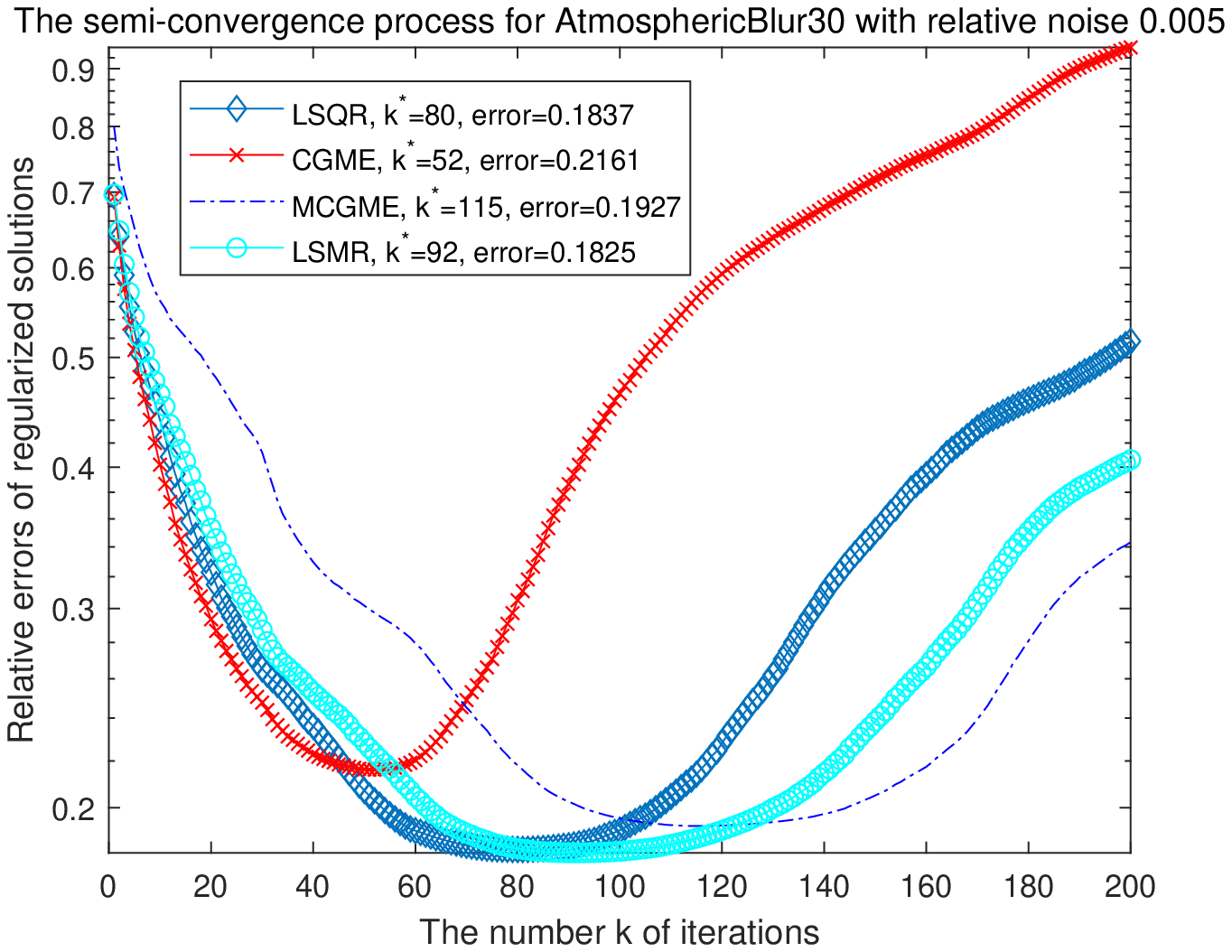}}
  \centerline{(b)}
\end{minipage}
\caption{{\rm (a)}: {\sf AtmosphericBlur10} and {\rm (b)}:
{\sf AtmosphericBlur30} with $\varepsilon=5\times
10^{-3}$.} \label{fig1}
\end{figure}

\begin{figure}
\begin{minipage}{0.48\linewidth}
  \centerline{\includegraphics[width=6.0cm,height=4.5cm]{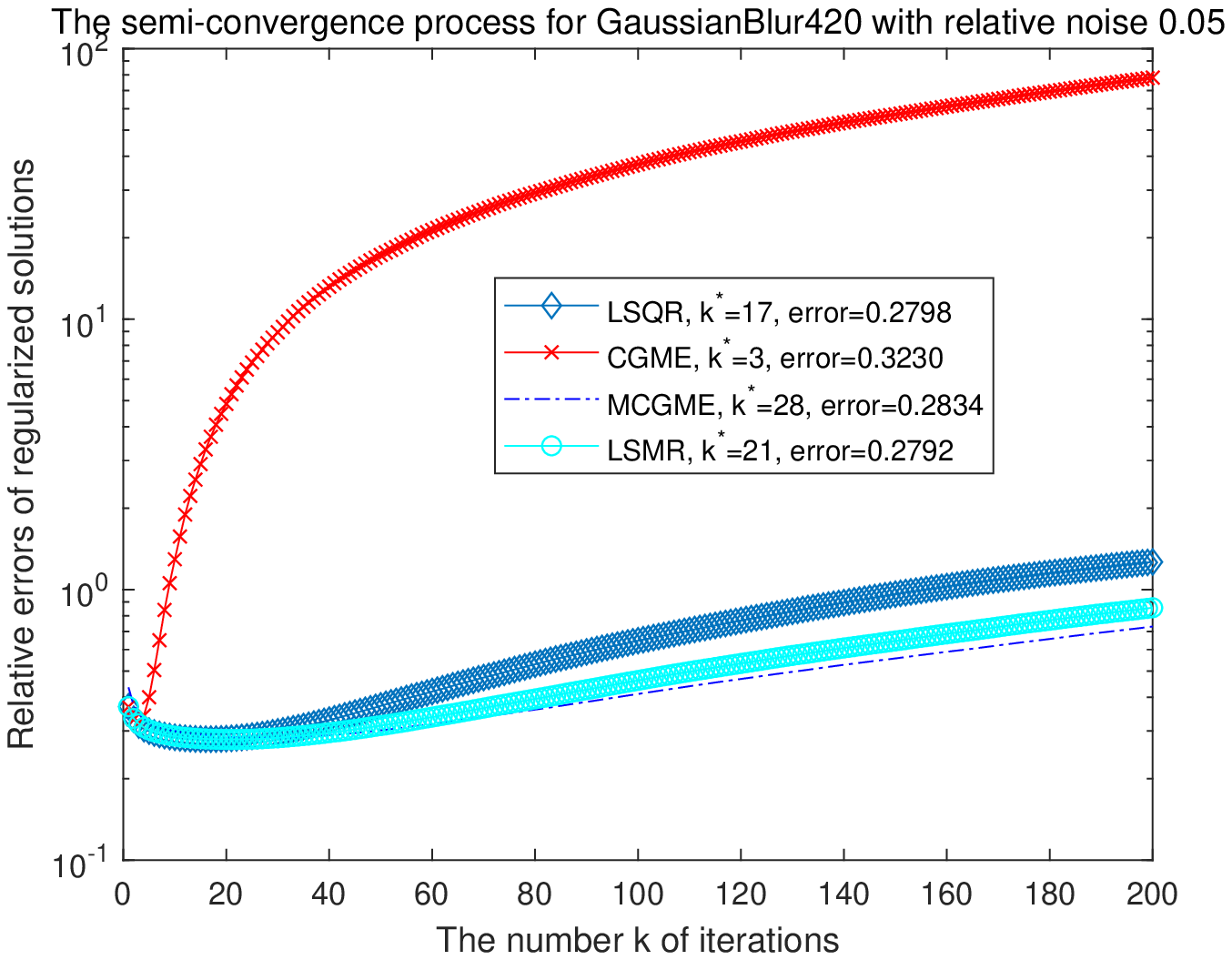}}
  \centerline{(a)}
\end{minipage}
\hfill
\begin{minipage}{0.48\linewidth}
  \centerline{\includegraphics[width=6.0cm,height=4.5cm]{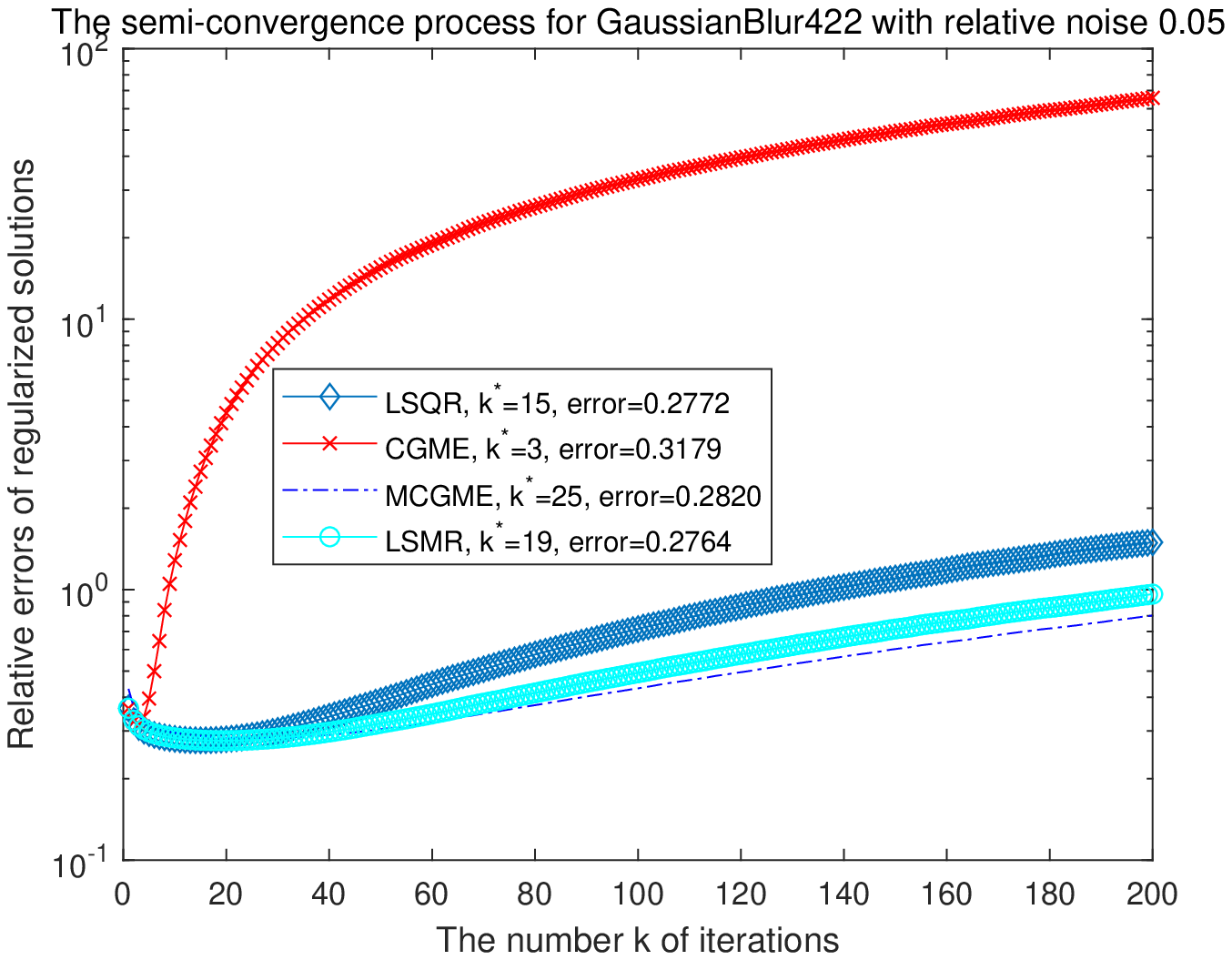}}
  \centerline{(b)}
\end{minipage}
\caption{{\rm (a)}: {\sf GaussianBlur420} and {\rm (b)}:
{\sf GaussianBlur422} with $\varepsilon=5\times
10^{-3}$}\label{fig2}
\end{figure}

\begin{figure}
\begin{minipage}{0.48\linewidth}
  \centerline{\includegraphics[width=6.0cm,height=4.5cm]{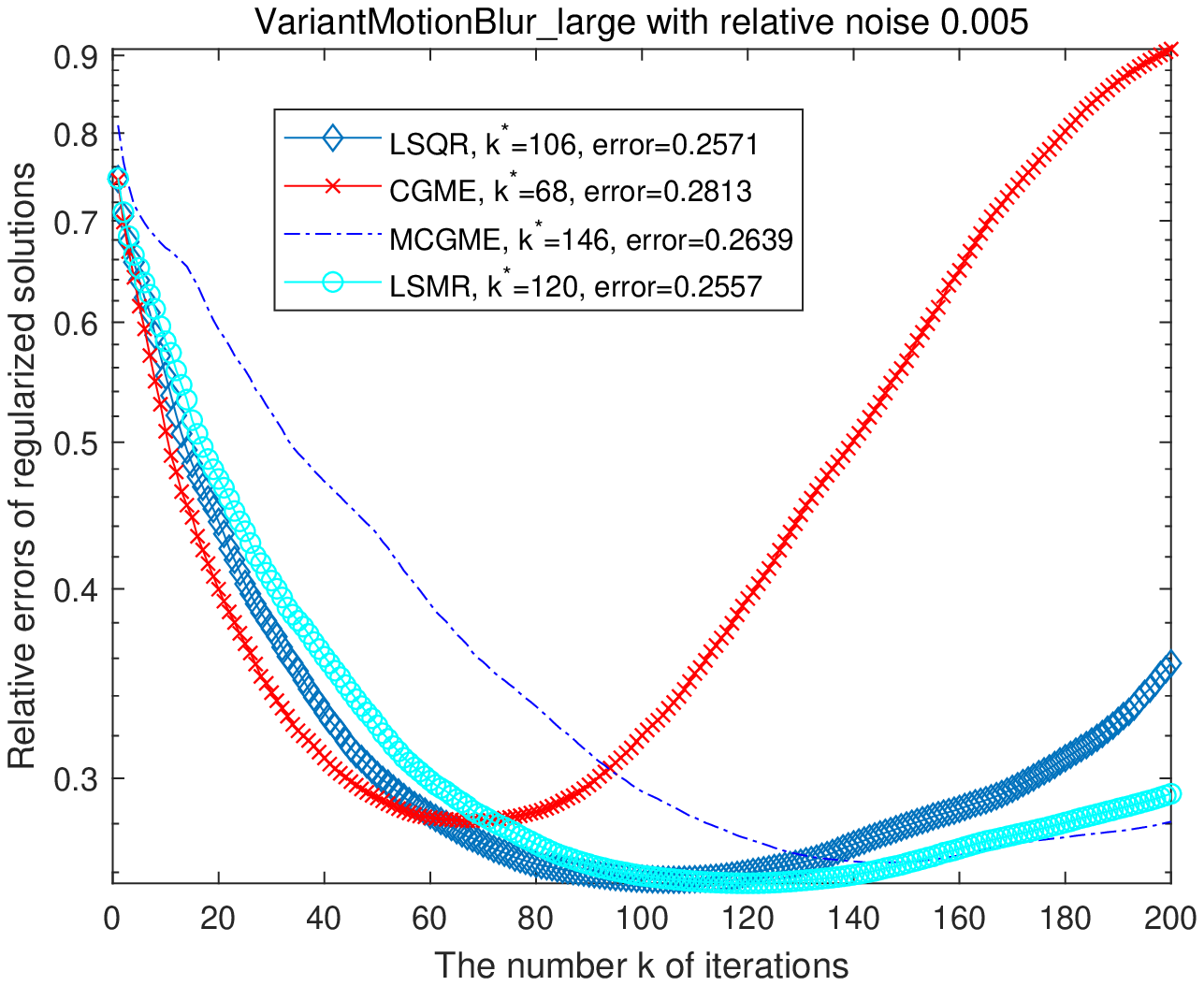}}
  \centerline{(a)}
\end{minipage}
\hfill
\begin{minipage}{0.48\linewidth}
  \centerline{\includegraphics[width=6.0cm,height=4.5cm]{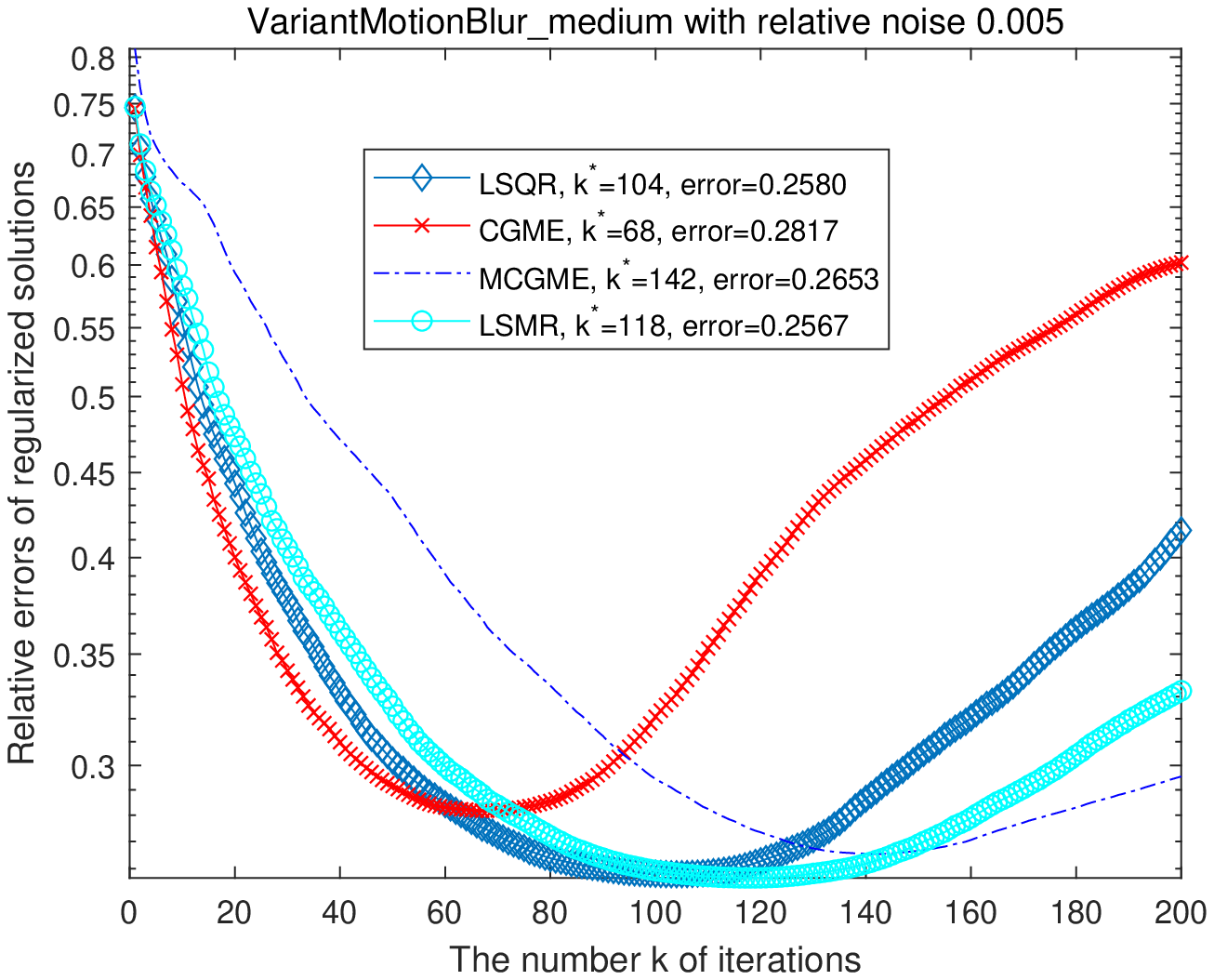}}
  \centerline{(b)}
\end{minipage}
\caption{{\rm (a)}: {\sf VariantMotionBlur\_large} and {\rm (b)}:
{\sf VariantMotionBlur\_medium} with $\varepsilon=5\times
10^{-3}$.} \label{fig3}
\end{figure}

\begin{figure}
\begin{minipage}{0.48\linewidth}
  \centerline{\includegraphics[width=6.0cm,height=4.5cm]{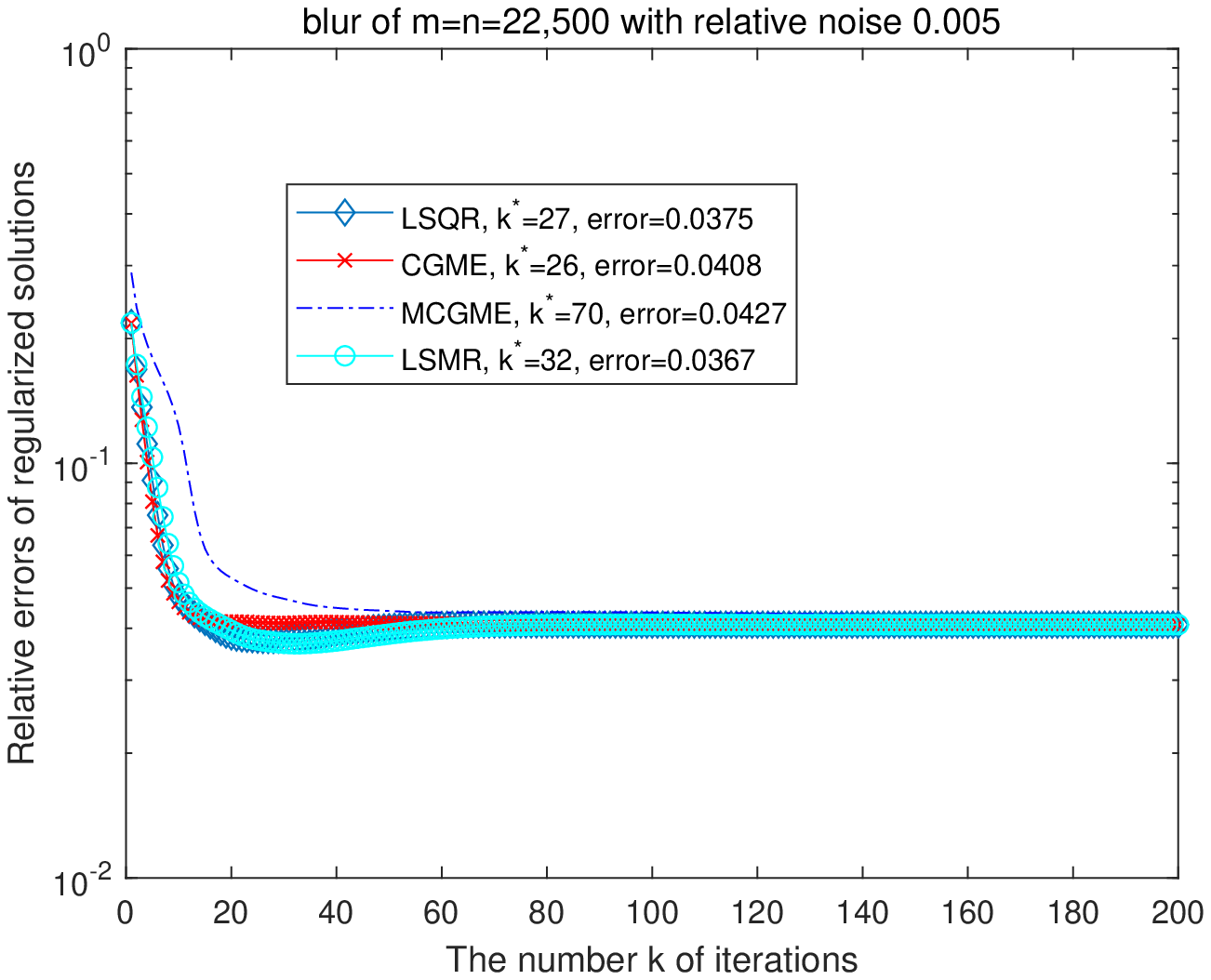}}
  \centerline{(a)}
\end{minipage}
\hfill
\begin{minipage}{0.48\linewidth}
  \centerline{\includegraphics[width=6.0cm,height=4.5cm]{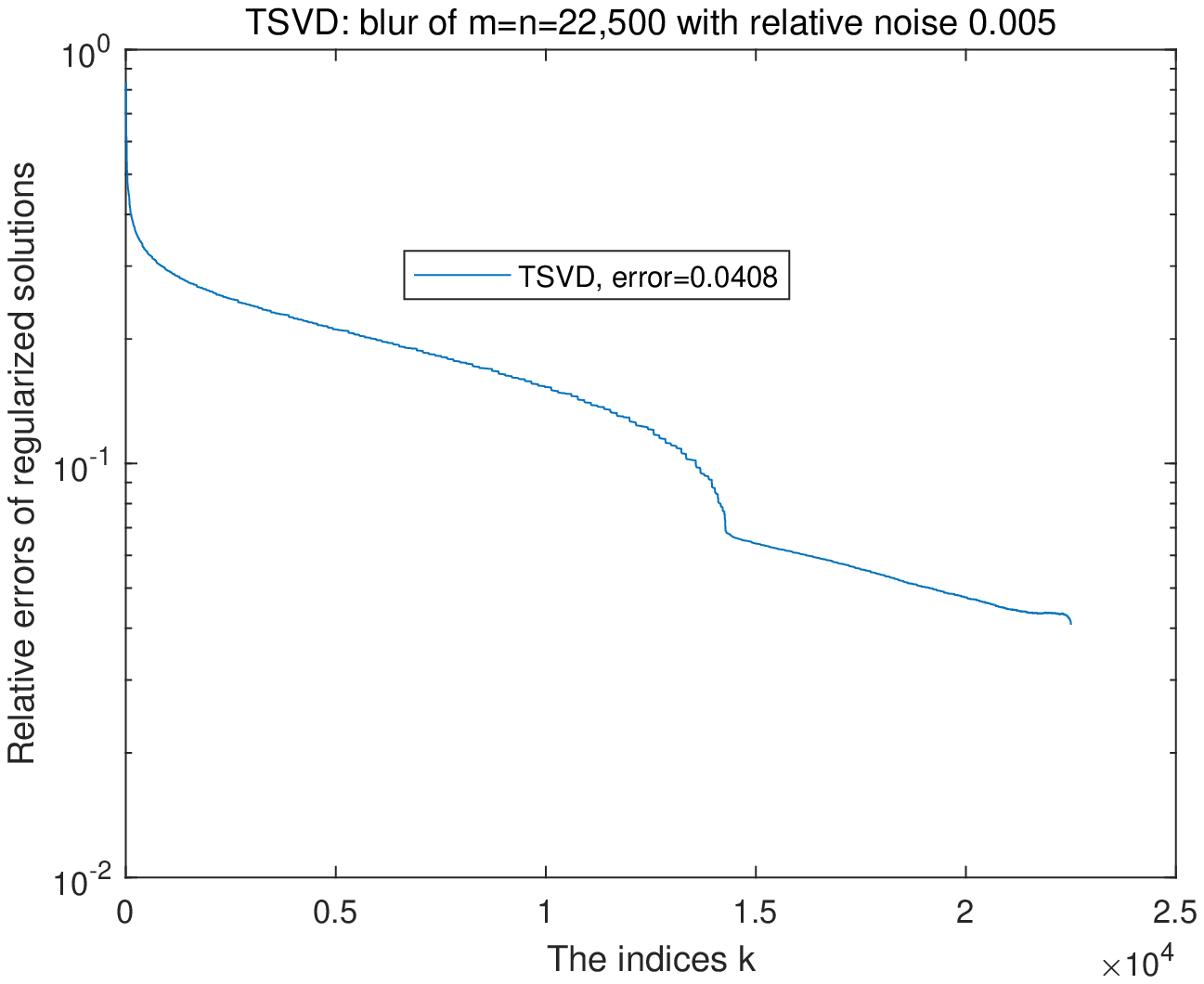}}
  \centerline{(b)}
\end{minipage}
\caption{{\sf blur}.} \label{fig4}
\end{figure}

\begin{figure}
\begin{minipage}{0.48\linewidth}
  \centerline{\includegraphics[width=6.0cm,height=4.5cm]{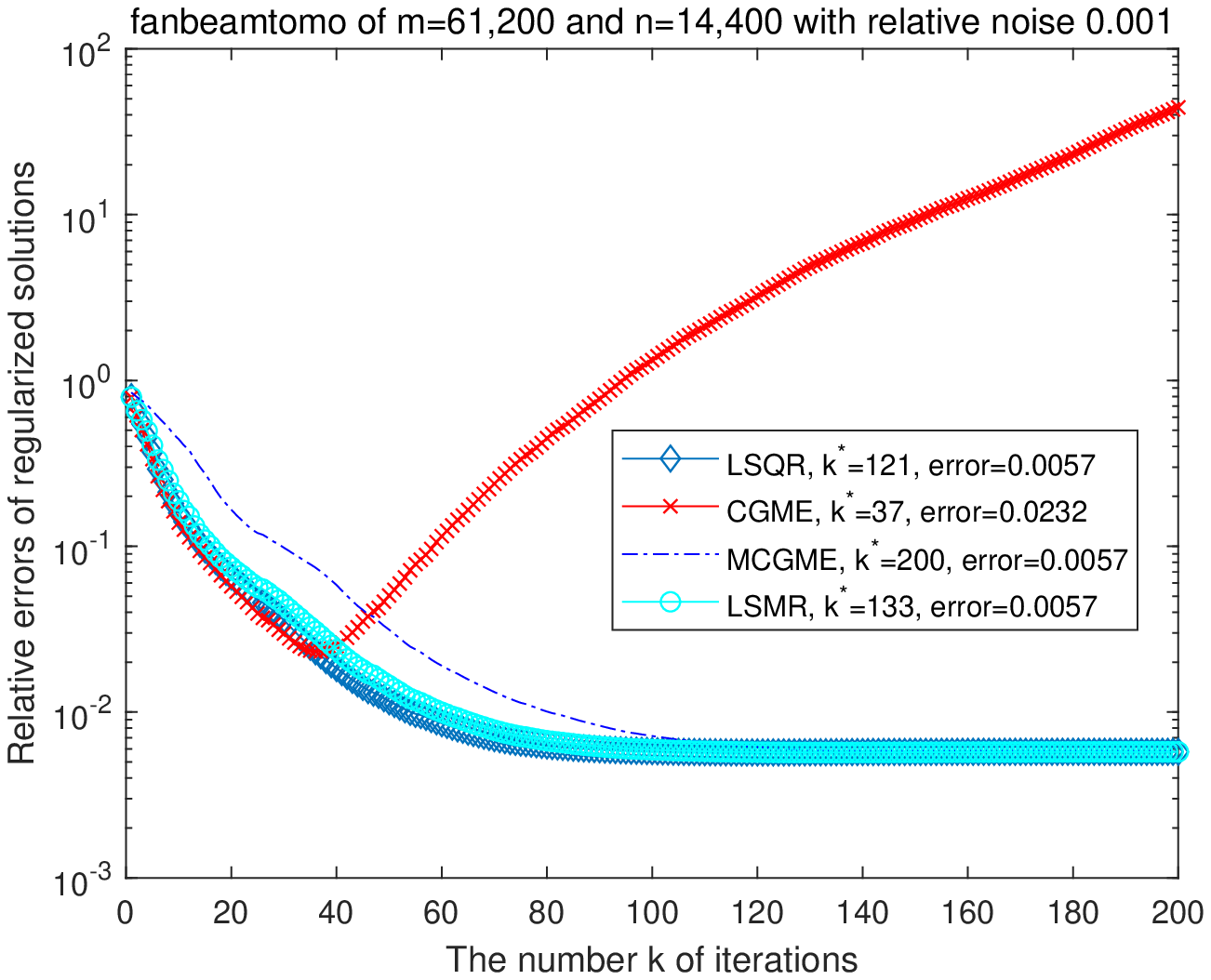}}
  \centerline{(a)}
\end{minipage}
\hfill
\begin{minipage}{0.48\linewidth}
  \centerline{\includegraphics[width=6.0cm,height=4.5cm]{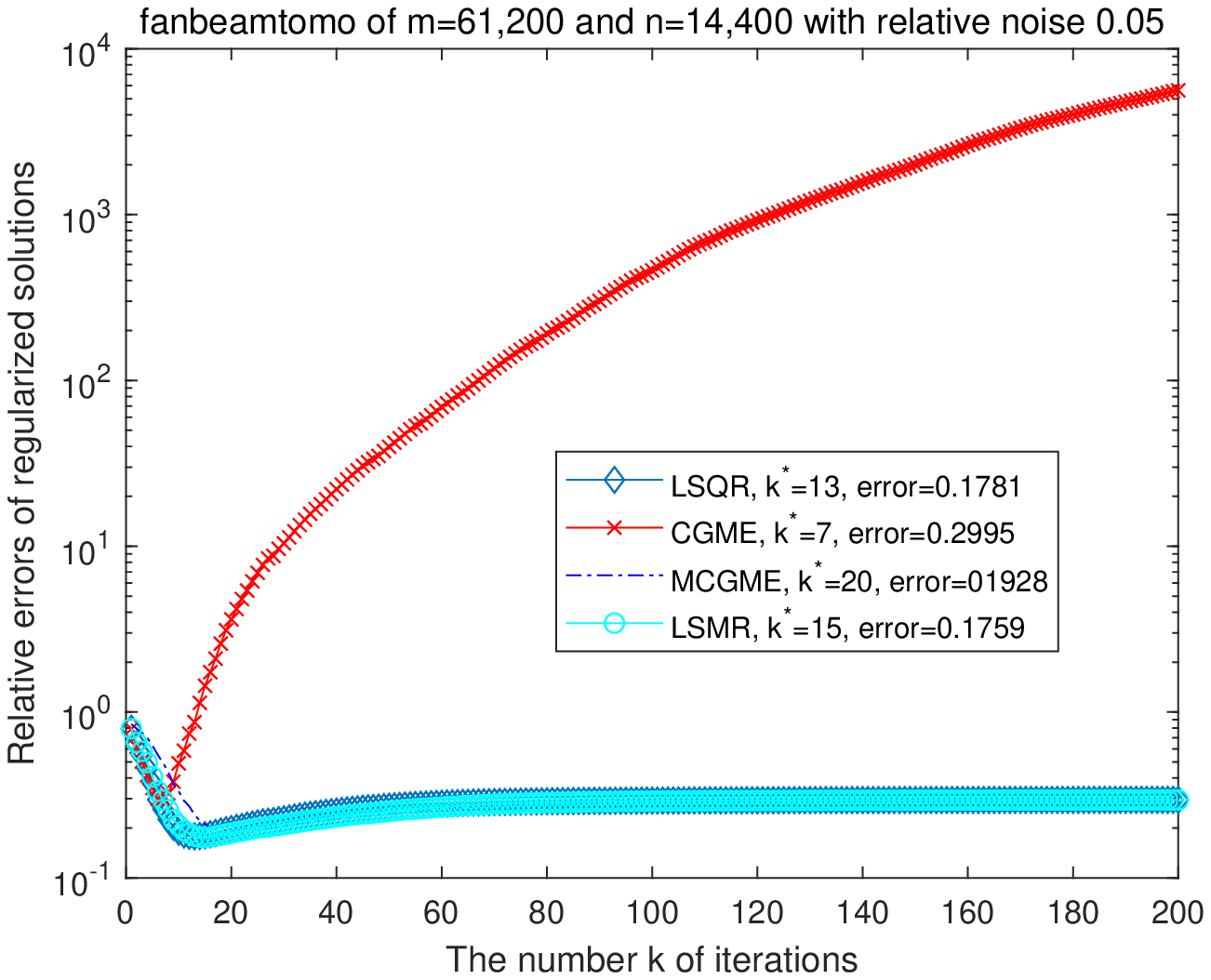}}
  \centerline{(b)}
\end{minipage}
\caption{{\sf fanbeamtomo} with $\varepsilon=10^{-3}$ and $5\times 10^{-2}$.} \label{fig5}
\end{figure}

\begin{figure}
\begin{minipage}{0.48\linewidth}
  \centerline{\includegraphics[width=6.0cm,height=4.5cm]{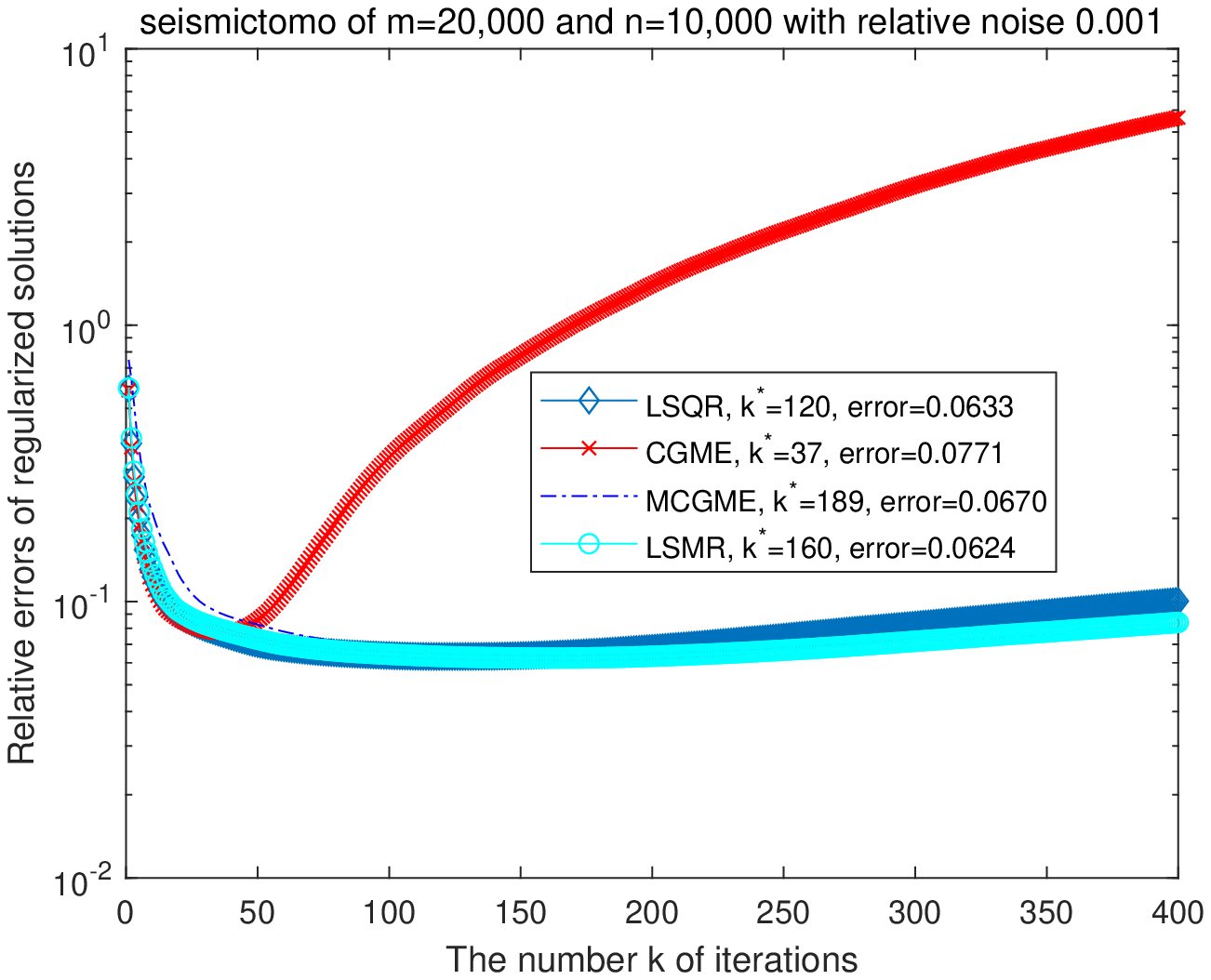}}
  \centerline{(a)}
\end{minipage}
\hfill
\begin{minipage}{0.48\linewidth}
  \centerline{\includegraphics[width=6.0cm,height=4.5cm]{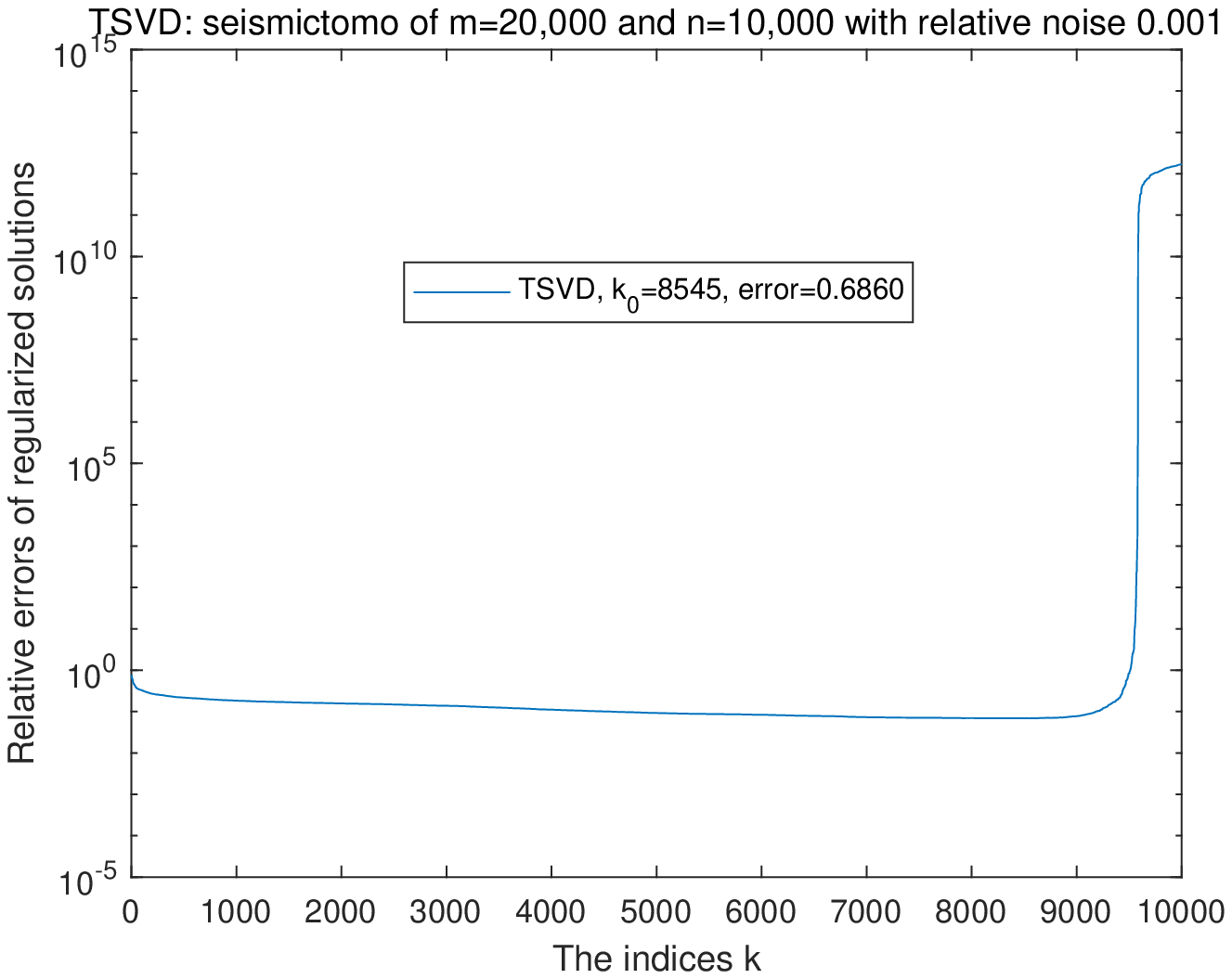}}
  \centerline{(b)}
\end{minipage}
\caption{{\sf seismictomo} with $\varepsilon=10^{-3}$.} \label{fig6}
\end{figure}

We now comment the figures and the related details in order.

Firstly, for all the problems in Table~\ref{tab1},
the semi-convergence of CGME occurs
earlier than LSQR and can be much earlier.
This confirms Theorem~\ref{semi}. The much earlier
semi-convergence of CGME indicates that
$\bar{\theta}_k^{(k)}<\sigma_{k_0+1}$ occurs much earlier
for CGME than $\theta_k^{(k)}<\sigma_{k_0+1}$ for LSQR.

Secondly, for all the problems, the best regularized solutions
$x_{k^*}^{cgme}$ are correspondingly less accurate than $x_{k^*}^{lsqr}$
considerably except for {\sf blur} in Figure~\ref{fig4}, where
the best regularized solution by CGME is almost as accurate as those
by LSQR, LSMR and MCGME. For all the 1D problems but {\sf baart}
and the 2D problem {\sf fanbeamtomo} with
$\varepsilon=10^{-3}$, the relative errors of the best regularized
solutions by CGME are twice to five times larger than the counterparts by
the other three ones, indicating that the regularization ability
is considerably inferior to the other three ones, given
that the relative errors by LSQR, LSMR and MCGME themselves
are only roughly $0.01\sim 0.1$;
see Figures~\ref{fig7} (a) and \ref{fig5} (a).
These results confirm Theorems~\ref{cgmeappr}--\ref{interlace}
and the analysis on them. We will make more comments on Figure~\ref{fig4} later.

Thirdly, for each of the problems, by a careful observation and comparison, we
have found that $x_{k}^{cgme}$ is more accurate than and at least at least as
accurate as $x_{k}^{lsqr}$ until the occurrence of CGME, after which
LSQR continues improving iterates until the occurrence of its
semi-convergence, as is clearly seen from
Figures~\ref{fig7}--\ref{fig6}.
These results justify our arguments on \eqref{accurcomp}.

Fourthly, for each of the 2D problems, the best regularized
solution $x_{k^*}^{lsmr}$ is at least as accurate as $x_{k^*}^{lsqr}$,
and the semi-convergence of LSMR always occurs no sooner and actually
later than that of LSQR. We notice that the relative error
of $x_{k^*}^{lsmr}$ is only slightly smaller
than that of $x_{k^*}^{lsqr}$, and there is
little difference between them. For all the 1D problems,
the semi-convergence of LSMR and LSQR occurs exactly at the same iterations,
and the best regularized solutions obtained by them have the same accuracy.
These results confirm Remark~\ref{semilsmr} and justify
that LSMR has the same regularization ability as that of LSQR.

Fifthly, for each of the test problems, MCGME improves CGME substantially. As a
matter of fact, for the 1D problems, the best regularized solutions by MCGME
have the same accuracy as those by LSQR and LSMR; for the 2D problems, the
best regularized solutions $x_{k^*}^{mcgme}$ are almost as accurate as
$x_{k^*}^{lsqr}$ and $x_{k^*}^{lsmr}$.

Sixthly, as we have stated, {\sf blur} and {\sf fanbeamtomo} are quite
well conditioned. With the relatively small $\varepsilon=10^{-3}$,
we observe from Figures~\ref{fig4}--\ref{fig5} that there is
no semi-convergence phenomenon for LSQR, LSMR and MCGME
as well as the TSVD method. This means that $e$
does not plays a part in regularization and
these methods solve these two problems as if they
were ordinary linear systems. Furthermore, it is clear from the figures
that the relative errors of regularized solutions obtained
by LSQR, LSMR and MCGME stabilize after 30 iterations for {\sf blur}
and 80 iterations for {\sf fanbeamtomo}, respectively.
Figures~\ref{fig4} (a) and \ref{fig5} (a) seems to indicate
that CGME has no semi-convergence phenomenon for the {\em square} {\sf blur}
and given $\varepsilon$ but it has for the
{\em rectangular} {\sf fanbeamtomo}.
However, this semi-convergence is in disguise and is not caused by the
noise $e$: For the rectangular {\sf fanbeamtomo},
\eqref{thetak2}, its proof and the analysis on it state that the smallest
singular value $\bar{\theta}_k^{(k)}$
of $\bar{B}_k$ can be arbitrarily small and approaches zero as $k$ increases.
As we have elaborated, $(\bar{\theta}_k^{(k)})^2$
approaches the eigenvalue zero of $AA^T$ as $k$ increases.
As a result, the projected problem
$\bar{B}_k y_k^{cgme}=\beta e_1^{(k)}$ involved in
CGME can become even worse conditioned than \eqref{eq1} itself
as $k$ increases for $A$ rectangular, causing that
$\|x_k^{cgme}\|$, which equals $\|Q_k y_k^{cgme}\|=\|y_k^{cgme}\|$,
and the relative error
$\frac{\|x_k^{cgme}-x_{true}\|}{\|x_{true}\|}$ tends to
infinity with respect to $k$. This can also be seen from \eqref{filter},
where we can easily check that $|f_k^{(k,cgme)}|\rightarrow \infty$
as $k$ increases since $\sigma_k$ is a constant
but $\bar{\theta}_k^{(k)}\rightarrow 0$ as $k$ increases.

In contrast,
the smallest singular values of the projection matrices are always
bounded from below by either $\sigma_n$ for LSQR (cf. \eqref{thetak1})
and MCGME (cf. \eqref{bkbarbk}) or $\sigma_n^2$
for LSMR (cf. \eqref{lsmrqr}), no matter how $A$ is rectangular or square.
This is why CGME has seemingly semi-convergence phenomenon for $A$
rectangular when the other solvers do not have. In the meantime, we see
that the best regularized solution by CGME is substantially less accurate
than those by the other three algorithms
for {\sf fanbeamtomo}. For the square {\sf blur} with $\varepsilon=10^{-3}$,
we see that the four Krylov solvers and the TSVD method
do not exhibit semi-convergence and compute
the solutions with very comparable accuracy. These results and
analysis tell us that CGME is definitely not a good choice
when $A$ is rectangular.

Seventhly, if the relative noise level $\varepsilon$ is increased
to $\varepsilon=0.05$,
the semi-convergence of LSQR, LSMR and MCGME occurs for {\sf fanbeamtomo},
as is seen from Figure~\ref{fig5}. We have also observed the semi-convergence
of the four algorithms and the TSVD method for {\sf blur} with $\varepsilon=0.05$.
We find that the best regularized solutions by LSQR, LSMR and MCGME have
very comparable accuracy but CGME computes a less accurate best regularized
solution.  We omit the corresponding figure. For the test problems,
we have also observed that the semi-convergence of
the TSVD method occurs much later than the four Krylov solvers,
i.e., $k^*\ll k_0$.

\section{Conclusions}\label{conclu}

For a general large-scale ill-posed problem \eqref{eq1}, iterative solvers
are only computationally viable. Of them, the Krylov solvers
LSQR, CGLS, CGME and LSMR have been commonly used.
In terms of the accuracy of the rank $k$ approximation to $A$ in LSQR,
in this paper we have derived accurate
estimates for the accuracy of the rank $k$ approximations to $A$ and
$A^TA$ that are involved in CGME and LSMR, respectively. We have made
detailed analyses on the approximation behavior of the singular
values of the projection matrices associated
with CGME and LSMR. In the meantime, we have
derived the filtered SVD expansion of CGME regularized iterates.
In conclusion, we have shown that the regularization of CGME is generally
inferior to LSQR and the semi-convergence of CGME occurs no later than
that of LSQR. We have extracted a best possible rank $k$ approximation to $A$
from the rank $(k+1)$ approximation $P_{k+1}P_{k+1}^TA$,
and have shown why such approximation is
as accurate as the rank $k$ approximation in LSQR.
Based on this analysis, as a by-product,
we have proposed a modified CGME (MCGME) method
that improves CGME substantially and
has the same regularization ability as LSQR.

We have substantially improved a fundamental result,
Theorem 9.3 in  \cite{halko11}, which gives a bound for the approximation accuracy
of the truncated rank $k$ SVD approximation to $A$ generated by randomized
algorithms and lacks a complete understanding to its considerable overestimate.
Our new bounds are unconditionally superior to theirs and reveal how the truncation
step affects the accuracy of the truncated rank $k$ approximation to $A$.

In the meantime, we have proved that LSMR has the same regularization
ability as LSQR and the semi-convergence of LSMR
occurs no sooner than that of LSQR. Particularly, we have shown that
LSMR has the full regularization for severely and moderately ill-posed problems
with suitable $\rho>1$ and $\alpha>1$.

We have made detailed numerical experiments to confirm our regularization results
on CGME and LSMR. We have also numerically demonstrated that the best regularized
solutions by MCGME are very comparable to those by LSQR.




\begin{thebibliography}{10}

\bibitem{aster}
R.~C. Aster, B.~Borchers, and C.~H. Thurber, \emph{Parameter {E}stimation and
  {I}nverse {P}roblems}, second ed., Elsevier, New York, 2013.

\bibitem{berisha}
S.~Berisha and J.~G. Nagy, \emph{Restore tools: Iterative methods for image
  restoration}, 2012, available from
  http://www.mathcs.emory.edu/$^\sim$nagy/RestoreTools.

\bibitem{bjorck96}
{\AA}.~Bj{\"{o}}rck, \emph{Numerical {M}ethods for {L}east {S}quares
  {P}roblems}, SIAM, Philadelphia, PA, 1996.

\bibitem{bjorck15}
\sameauthor, \emph{Numerical {M}ethods in {M}atrix {C}omputations}, Texts in
  Applied Mathematics, vol.~59, Springer, Cham, 2015.

\bibitem{chung15}
J.~Chung and K.~Palmer, \emph{A hybrid {LSMR} algorithm for large-scale
  {T}ikhonov regularization}, SIAM J. Sci. Comput., 37 (5) (2015),
  pp.~S562--S580.

\bibitem{craig}
E.~J. Craig, \emph{The {$N$}-step iteration procedures}, J. Math. Phys.
  34 (1955), pp.~64--73.

\bibitem{demmel}
J.~Demmel, \emph{Applied {N}umerical {L}inear {A}lgebra}, SIAM,
  Philadelphia, PA, 1997.

\bibitem{eicke}
B.~Eicke, A.~K. Lious, and R.~Plato, \emph{The instability of some gradient
  methods for ill-posed problems}, Numer. Math., 58 (1) (1990), pp.~129--134.

\bibitem{engl93}
H.~W. Engl, \emph{Regularization methods for the stable solution of inverse
  problems}, Surveys Math. Indust., 3 (2) (1993), pp.~71--143.

\bibitem{engl00}
H.~W. Engl, M.~Hanke, and A.~Neubauer, \emph{{R}egularization of {I}nverse
  {P}roblems}, Kluwer Academic Publishers, 2000.

\bibitem{firro97}
R.~D. Fierro, G.~H. Golub, P.~C. Hansen, and D.~P. O'Leary,
  \emph{Regularization by truncated total least squares}, SIAM J. Sci. Comput.,
  18 (4) (1997), pp.~1223--1241.

\bibitem{fong}
D.~C.~L. Fong and M.~Saunders, \emph{L{SMR}: an iterative algorithm for sparse
  least-squares problems}, SIAM J. Sci. Comput., 33 (5) (2011),
  pp.~2950--2971.

\bibitem{gazzola15}
S.~Gazzola and P.~Novati, \emph{Inheritance of the discrete {P}icard condition
  in {K}rylov subspace methods}, BIT Numer. Math., 56 (3) (2016), pp.~893--918.

\bibitem{gilyazov}
S.~F. Gilyazov and N.~L. Gol'dman, \emph{Regularization of {I}ll-{P}osed
  {P}roblems by {I}teration {M}ethods}, Kluwer Academic Publishers, Dordrecht, 2000.

\bibitem{golub89}
G.~H. Golub and D.~P. O'Leary, \emph{Some history of the conjugate gradient and
  {L}anczos algorithms: 1948--1976}, SIAM Rev., 31 (1) (1989),
  pp.~50--102.

\bibitem{halko11}
N.~Halko, P.~G. Martinsson, and J.~A. Tropp, \emph{Finding structure with
  randomness: probabilistic algorithms for constructing approximate matrix
  decompositions}, SIAM Rev., 53 (2) (2011), pp.~217--288.

\bibitem{hanke95}
M.~Hanke, \emph{Conjugate {gr}adient {T}ype {M}ethods for {I}ll-{P}osed
  {P}roblems}, Pitman Research Notes in Mathematics Series, vol. 327, Longman,
  Essex, 1995.

\bibitem{hanke01}
\sameauthor, \emph{On {L}anczos based methods for the regularization of discrete
  ill-posed problems}, BIT Numer. Math., 41 (5) (2001), Suppl.,
  pp.~1008--1018.

\bibitem{hanke93}
M.~Hanke and P.~C. Hansen, \emph{Regularization methods for large-scale
  problems}, Surveys Math. Indust., 3 (4) (1993), pp.~253--315.

\bibitem{hansen90}
P.~C. Hansen, \emph{The discrete {P}icard condition for discrete ill-posed
  problems}, BIT, 30 (4) (1990), pp.~658--672.

\bibitem{hansen90b}
\sameauthor, \emph{Truncated singular value decomposition solutions to discrete
  ill-posed problems with ill-determined numerical rank}, SIAM J. Sci. Statist.
  Comput., 11 (3) (1990), pp.~503--518.

\bibitem{hansen98}
\sameauthor, \emph{Rank-{D}eficient and {D}iscrete {I}ll-{P}osed {P}roblems:
  {N}umerical {A}spects of {L}inear {I}nversion}, SIAM, Philadelphia, PA, 1998.


\bibitem{hansen07}
\sameauthor, \emph{Regularization {T}ools version 4.0 for {M}atlab 7.3}, Numer.
  Algor., 46 (2) (2007), pp.~189--194.

\bibitem{hansen10}
\sameauthor, \emph{Discrete {I}nverse {P}roblems: {I}nsight and {A}lgorithms},
  SIAM, Philadelphia, PA, 2010.

\bibitem{hansen12}
P.~C. Hansen and M.~Saxild-Hansen, \emph{{AIR} tools--a {MATLAB} package of
  algebraic iterative reconstruction methods}, J. Comput. Appl. Math.,
  236 (2012), pp.~2167--2178.

\bibitem{hestenes}
M.~R. Hestenes and E.~Stiefel, \emph{Methods of conjugate gradients for solving
  linear systems}, J. Res. Nat. Bur. Stand., 49 (1952), pp.~409--436.

\bibitem{hps16}
M.~R. Hn\v{e}tynkov\'{a}, Marie Kub\'{i}nov\'{a}, and M.~Ple\v{s}inger,
  \emph{Noise representation in residuals of {LSQR}, {LSMR}, and {C}raig
  regularization}, Linear Algebra Appl., 533 (2017), pp.~357--379.

\bibitem{hps09}
M.~R. Hn\v{e}tynkov\'{a}, M.~Ple\v{s}inger, and Z.~Strako\v{s}, \emph{The
  regularizing effect of the {G}olub-{K}ahan iterative bidiagonalization and
  revealing the noise level in the data}, BIT Numer. Math., 49 (4) (2009),
  pp.~669--696.

\bibitem{hofmann86}
B.~Hofmann, \emph{{R}egularization for {A}pplied {I}nverse and {I}ll-{P}osed
  {P}roblems}, Teubner, Stuttgart, Germany, 1986.

\bibitem{huangjia}
Y.~Huang and Z.~Jia, \emph{Some results on the regularization of {LSQR} for large-scale
  ill-posed problems}, Science China Math., 60 (4) (2017), pp.~701--718.


\bibitem{huangjia17}
\sameauthor, \emph{On regularizing effects of {MINRES} and {MR-II} for
  large-scale symmetric discrete ill-posed problems}, J. Comput. Appl. Math.,
  320 (2017), pp.~145--163.


\bibitem{jia18a}
Z.~Jia, \emph{{A}pproximation accuracy of the {K}rylov subspaces for linear
  discrete ill-posed problems},  (2018), arXiv:math.NA/1805.10132.

\bibitem{jia18b}
\sameauthor, \emph{The low rank approximations and {R}itz values in {LSQR} for
  linear discrete ill-posed problems},  (2018), arXiv:math.NA/1811.03454.


\bibitem{jia03}
Z.~Jia and D.~Niu, \emph{An implicitly restarted refined bidiagonalization
  {L}anczos method for computing a partial singular value decomposition}, SIAM
  J. Matrix Anal. Appl., 25 (1) (2003), pp.~246--265.

\bibitem{kaipio}
J.~Kaipio and E.~Somersalo, \emph{Statistical and {C}omputational {I}nverse
  {P}roblems}, Springer-Verlag, New
  York, 2005.

\bibitem{kern}
M.~Kern, \emph{Numerical {M}ethods for {I}nverse {P}roblems}, John Wiley \&
  Sons, Inc., 2016.

\bibitem{kirsch}
A.~Kirsch, \emph{An {I}ntroduction to the {M}athematical {T}heory of {I}nverse
  {P}roblems}, second ed., Springer, New
  York, 2011.

\bibitem{meurant}
G.~Meurant, \emph{The {L}anczos and {C}onjugate {G}radient {A}lgorithms: From
  {T}heory to {F}inite {P}recision {C}omputations}, SIAM,
  Philadelphia, PA, 2006.

\bibitem{natterer}
F.~Natterer, \emph{The {M}athematics of {C}omputerized {T}omography},
SIAM, Philadelphia, PA, 2001.

\bibitem{paige75}
C.~C. Paige and M.~A. Saunders, \emph{Solutions of sparse indefinite systems of
  linear equations}, SIAM J. Numer. Anal., 12 (4) (1975), pp.~617--629.

\bibitem{paige82}
\sameauthor, \emph{{LSQR}: an algorithm for sparse linear equations and sparse
  least squares}, ACM Trans. Math. Software, 8 (1) (1982), pp.~43--71.

\bibitem{paige06}
C.~C. Paige and Z.~Z. Strako\v{s}, \emph{Core problems in linear algebraic
  systems}, SIAM J. Matrix Anal. Appl., 27 (3) (2005), pp.~861--875.

\bibitem{parlett}
B.~N. Parlett, \emph{The {S}ymmetric {E}igenvalue {P}roblem},
SIAM, Philadelphia, PA, 1998.

\bibitem{stewartsun}
G.~W. Stewart and J.-G Sun, \emph{Matrix {P}erturbation {T}heory},
Academic Press, Inc., Boston, MA, 1990.


\bibitem{vorst90}
A.~van~der Sluis and H.~A. van~der Vorst, \emph{S{IRT}- and {CG}-type methods
  for the iterative solution of sparse linear least-squares problems}, Linear
  Algebra Appl., 130 (1990), pp.~257--303.

\bibitem{varah79}
J.~M. Varah, \emph{A practical examination of some numerical methods for linear
  discrete ill-posed problems}, SIAM Rev., 21 (1) (1979), pp.~100--111.

\bibitem{vogel02}
C.~R. Vogel, \emph{Computational {M}ethods for {I}nverse {P}roblems},
SIAM, Philadelphia, PA, 2002.

\bibitem{wilkinson}
J.~H. Wilkinson, \emph{The {A}lgebraic {E}igenvalue {P}roblem}, Clarendon
  Press, Oxford, 1965.

\end{thebibliography}

\end{document}